\documentclass[a4paper, 11pt, twoside]{article}

\usepackage{amsmath, amscd, amsfonts, amssymb, amsthm, latexsym, url, color, todonotes, bm, framed}
\usepackage{graphicx}
\usepackage[left=1.2in, right=1.2in, top=1.3in, bottom=1.3in, includefoot, headheight=13.6pt]{geometry}
\usepackage{mathtools}
\input{xy}
\xyoption{all}

\usepackage[ps2pdf=true,colorlinks,breaklinks=true]{hyperref}
\usepackage[figure,table]{hypcap}
\hypersetup{
	bookmarksnumbered,
	pdfstartview={FitH},
	citecolor={black},
	linkcolor={black},
	urlcolor={black},
	pdfpagemode={UseOutlines}
}
\makeatletter
\newcommand\org@hypertarget{}
\let\org@hypertarget\hypertarget
\renewcommand\hypertarget[2]{%
  \Hy@raisedlink{\org@hypertarget{#1}{}}#2%
} 
\makeatother 

\newtheorem{theorem}{Theorem}[section]
\newtheorem{lemma}[theorem]{Lemma}
\newtheorem{corollary}[theorem]{Corollary}
\newtheorem{proposition}[theorem]{Proposition}

\theoremstyle{definition}

\newtheorem{remark}[theorem]{Remark}
\newtheorem{example}[theorem]{Example}
\newtheorem{conjecture}[theorem]{Conjecture}

\newcommand{\xysquare}[8]{
\[\xymatrix{
#1 \ar@{#5}[r] \ar@{#6}[d] & #2 \ar@{#7}[d]\\
#3 \ar@{#8}[r] & #4
}\]
}

\DeclareMathOperator*{\holim}{\operatorname*{holim}}

\newcommand{\al}{\alpha}
\newcommand{\bb}{\mathbb}

\newcommand{\blob}{\bullet}

\newcommand{\comment}[1]{}

\newcommand{\ep}{\varepsilon}

\newcommand{\into}{\hookrightarrow}
\newcommand{\isoto}{\stackrel{\simeq}{\to}}

\newcommand{\onto}{\twoheadrightarrow}
\newcommand{\op}{\operatorname}
\newcommand{\pid}[1]{\langle #1 \rangle}

\newcommand{\res}{\overline}
\newcommand{\roi}{\mathcal{O}}

\newcommand{\sub}[1]{{\mbox{\scriptsize #1}}}

\newcommand{\To}{\longrightarrow}

\newcommand{\xto}{\xrightarrow}

\newcommand{\TC}{T\!C}
\newcommand{\NS}{N\!S}

\renewcommand{\cal}{\mathcal}
\renewcommand{\hat}{\widehat}

\newcommand{\indlim}{\varinjlim}
\renewcommand{\tilde}{\widetilde}
\renewcommand{\Im}{\operatorname{Im}}
\renewcommand{\ker}{\operatorname{Ker}}
\renewcommand{\projlim}{\varprojlim}

\DeclareMathOperator{\Frac}{Frac}
\DeclareMathOperator{\Gal}{Gal}
\DeclareMathOperator{\Hom}{Hom}

\DeclareMathOperator{\Pic}{Pic}

\DeclareMathOperator{\Spec}{Spec}

\newcommand{\CH}{C\!H}

\usepackage{fancyhdr}

\usepackage{sectsty}
\sectionfont{\Large\sc\centering}
\chapterfont{\large\sc\centering}
\chaptertitlefont{\LARGE\centering}
\partfont{\centering}

\begin{document}

\title{A Variational Tate Conjecture in crystalline cohomology}

\author{Matthew Morrow}

\date{}

\maketitle


\begin{abstract}
Given a smooth, proper family of varieties in characteristic $p>0$, and a cycle $z$ on a fibre of the family, we consider a Variational Tate Conjecture characterising, in terms of the crystalline cycle class of $z$, whether $z$ extends cohomologically to the entire family. This is a characteristic $p$ analogue of Grothendieck's Variational Hodge Conjecture. We prove the conjecture for divisors, and an infinitesimal variant of the conjecture for cycles of higher codimension.

This can be used to reduce the $\ell$-adic Tate conjecture for divisors over finite fields to the case of surfaces.
\end{abstract}

\tableofcontents

\setcounter{section}{-1}
\section{Introduction and statement of main results}
Let $f:X\to S$ be a smooth, proper morphism of smooth varieties over a field $k$, and let $s\in S$ be a closed point. Grothendieck's Variational Hodge or $\ell$-adic Tate Conjecture \cite[pg.~359]{Grothendieck1966} gives conditions on the cohomology class of an algebraic cycle on $X_s$ under which the cycle conjecturally extends cohomologically to the entire family $X$. According as $k$ has characteristic $0$ or $p>0$, the cohomology theory used to formulate Grothendieck's conjecture is de Rham or $\ell$-adic \'etale ($\ell\neq p$). The Variational Hodge Conjecture is easily proved for divisors using the exponential map; in contrast, it seems little is known about the $\ell$-adic Variational Tate Conjecture.

Suspecting that deformation problems in characteristic $p$ are best understood $p$-adically, we study in this article a Variational Tate Conjecture in crystalline cohomology. To state it, we now fix some notation for the rest of the article. Unless stated otherwise, $k$ is a perfect field of characteristic $p>0$, and $W=W(k)$, $K=\Frac W$. Given any reasonable (not necessarily of finite type) scheme $X$ over $k$ of characteristic $p$, let $H^n_\sub{crys}(X)=H^n_\sub{crys}(X/W)\otimes_WK$ denote its rational crystalline cohomology. Assuming now that $X$ is a smooth $k$-variety, let $cl_i:\CH^i(X)_\bb Q\To H^{2i}_\sub{crys}(X)$ denote the crystalline cycle class map for any $i\ge0$; these land in the $p^i$-eigenspace $H^{2i}_\sub{crys}(X)^{\phi=p^i}$ of the absolute Frobenius $\phi:x\mapsto x^p$.

The Variational Tate Conjecture in crystalline cohomology to be studied is as follows:

\begin{conjecture}[Crystalline Variational Tate Conjecture]\label{VTC}
Let $f:X\to S$ be a smooth, proper morphism of smooth $k$-varieties, $s\in S$ a closed point, and $z\in \CH^i(X_s)_\bb Q$. Let $c:=cl_i(z)\in H^{2i}_\sub{crys}(X_s)$. Then the following are equivalent:
\begin{description}
\item[(deform)] There exists $\tilde z\in \CH^i(X)_\bb Q$ such that $cl_i(\tilde z)|_{X_s}=c$.
\item[(crys)] $c$ lifts to $H^{2i}_\sub{crys}(X)$.
\item[(crys$\boldsymbol{-\phi}$)] $c$ lifts to $H^{2i}_\sub{crys}(X)^{\phi=p^i}$.
\item[(flat)] $c$ is flat, i.e., it lifts to $H^0_\sub{crys}(S,R^{2i}f_*\roi_{X/K})$.
\end{description}
\end{conjecture}

A more detailed discussion of Conjecture \ref{VTC}, including an equivalent formulation via rigid cohomology and an explanation of the condition (flat), is given in Section \ref{section_general_comments}.

Our first main result is the proof of Conjecture \ref{VTC} for divisors, at least assuming that the family is projective:

\begin{theorem}[See Thm.~\ref{theorem_VTC_for_divisors}]\label{theorem_2}
Conjecture \ref{VTC} is true for divisors, i.e.,~when $i=1$, if $f$ is projective.
\end{theorem}

A standard pencil argument yields a variant of Theorem \ref{theorem_2} for extending divisors from a hyperplane section of a projective variety; see Corollary \ref{corollary_hypersurfaces}. This is used to establish the following application, which is already known to certain experts in the field including A.~J.~de Jong \cite{deJong20??}, but unfortunately seems not to be more widely known:

\begin{corollary}[See Thm.~\ref{theorem_reduction_of_Tate}]\label{intro_corol_tate}
Assume that the Tate conjecture for divisors is true for all smooth, projective surfaces over a finite field $k$. Then the Tate conjecture for divisors is true for all smooth, projective varieties over $k$.
\end{corollary}

Regarding Corollary \ref{intro_corol_tate}, we remark that the ``Tate conjecture for divisors'' is independent of the chosen Weil cohomology theory (see Proposition \ref{proposition_Tate_conj}); in particular, the corollary may be stated in terms of $\ell$-adic \'etale cohomology, though the proof is implicitly a combination of crystalline and $\ell$-adic techniques.

We also prove a variant of Theorem \ref{theorem_2} for line bundles on smooth, projective schemes over the spectrum of a power series ring $k[[t_1,\dots,t_m]]$; see Theorem \ref{theorem_line_bundles}. Combining this with N.~Katz' results on slope filtrations of $F$-crystals over $k[[t]]$ yields the following consequence:

\begin{corollary}[See Thm.~\ref{theorem_constant_families}]\label{corollary_2a}
Let $X$ be a smooth, proper scheme over $k[[t]]$, where $k$ is an algebraically closed field of characteristic $p$. Assume that $R^nf_*\roi_{X/W}$ is locally-free for all $n\ge0$ and is a constant $F$-crystal for $n=2$. Then the cokernel of the restriction map $\Pic(X)\to\Pic(X\times_{k[[t]]}k)$ is killed by a power of $p$.
\end{corollary}

Corollary \ref{corollary_2a} applies in particular to supersingular families of K3 surfaces over $k[[t]]$, thereby reproving a result of M.~Artin \cite{Artin1974}.

We now explain our other variational results while simultaneously indicating the main ideas of the proofs. Firstly, Section \ref{section_partie_fixe} is devoted to the proof of a crystalline analogue of Deligne's Th\'eor\`eme de la Partie Fixe, stating that the crystalline Leray spectral sequence for the morphism $f$ degenerates in a strong sense; see Theorem \ref{theorem_crystalline_partie_fixe}. This implies, in the situation of Conjecture \ref{VTC}, that conditions (crys), (crys-$\phi$), and (flat) are in fact equivalent (assuming $f$ is projective); these conditions are also equivalent to their analogues in rigid cohomology, at least conditionally under expected hypotheses (see Remarks  \ref{remark_rigid_1} and \ref{remark_rigid_2}).

Hence, to prove Theorem \ref{theorem_2}, it is enough to show that (crys-$\phi$) implies (deform) for divisors. By standard arguments, we may base change by $k^\sub{alg}$, replace $S$ by $\Spec\hat\roi_{S,s}$ and, identifying divisors with line bundles, then prove the following, in which it is only necessary to invert $p$:

\begin{theorem}[{See Corol.~\ref{corollary_line_bundles}}]\label{theorem_3}
Let $X$ be a smooth, proper scheme over $A=k[[t_1,\dots,t_m]]$, where $k$ is an algebraically closed field of characteristic $p$, and let $L\in\Pic(X\times_Ak)[\tfrac1p]$. Then the following are equivalent:
\begin{enumerate}
\item There exists $\tilde{L}\in\Pic(X)[\tfrac1p]$ such that $\tilde L|_{X\times_Ak}=L$.
\item The first crystalline Chern class $c_1(L)\in H^2_\sub{crys}(X\times_Ak)$ lifts to $H^2_\sub{crys}(X)^{\phi=p}$.
\end{enumerate}
\end{theorem}

Theorem \ref{theorem_3} is an exact analogue in equal characteristic $p$ of the following existing deformation theorems in mixed characteristic and in equal characteristic zero:
\begin{itemize}
\item $p$-adic Variational Hodge Conjecture for line bundles (P.~Berthelot and A.~Ogus \cite[Thm.~3.8]{BerthelotOgus1983}).
Let $X$ be a smooth, proper scheme over a complete discrete valuation ring $V$ of mixed characteristic with perfect residue field $k$, and let $L\in\Pic(X\times_Vk)[\tfrac1p]$. Then $L$ lifts to $\Pic(X)[\tfrac1p]$ if and only if its first crystalline Chern class $c_1(L)\in H^2_\sub{crys}(X\times_Vk)$ belongs to $F^1H^2_\sub{dR}(X/V)[\tfrac1p]$ under the comparison isomorphism $H^2_\sub{crys}(X\times_Vk)\cong H^2_\sub{dR}(X/V)[\tfrac1p]$.
\item Local case of the Variational Hodge Conjecture for line bundles (folklore).
Let $X$ be a smooth, proper scheme over $A=k[[t_1,\dots,t_m]]$, where $k$ is a field of characteristic $0$, and let $L\in\Pic(X\times_Ak)$. Then $L$ lifts to $\Pic(X)$ if and only if its first de Rham Chern class $c_1(L)\in H^2_\sub{crys}(X\times_Ak)$ lifts to $F^1H^2_\sub{dR}(X/k)$, which denotes the Hodge filtration on the continuous de Rham cohomology of $X$ over~$k$ (alternatively, $F^1H^2_\sub{dR}(X/k)$ identifies with the subspace of $F^1H^2_\sub{dR}(X/A)$ annihilated by the Gauss--Manin connection).
\end{itemize}

We finally state our infinitesimal version of Conjecture \ref{VTC} which holds in all codimensions; it provides a necessary and sufficient condition under which the $K_0$ class of a vector bundle on the special fibre admits infinitesimal extensions of all orders:

\begin{theorem}[See Thm.~\ref{main_infinitesimal_theorem}]\label{theorem_4}
Let $X$ be a smooth, proper scheme over $A=k[[t_1,\dots,t_m]]$, where $k$ is a finite or algebraically closed field of characteristic $p$; let $Y$ denote the special fibre and $Y_r:=X\times_AA/\pid{t_1,\dots,t_m}^r$ its infinitesimal thickenings. Then the following are equivalent for any $z\in K_0(Y)[\tfrac1p]$:
\begin{enumerate}\itemsep0pt
\item $z$ lifts to $(\projlim_r K_0(Y_r))[\tfrac1p]$.
\item The crystalline Chern character $ch(z)\in\bigoplus_{i\ge0}H_\sub{crys}^{2i}(Y)$ lifts to $\bigoplus_{i\ge0}H^{2i}_\sub{crys}(X)^{\phi=p^i}$.
\end{enumerate}
\end{theorem}

Theorem \ref{theorem_4} is an analogue in characteristic $p$ of a deformation result in mixed characteristic due to S.~Bloch, H.~Esnault, and M.~Kerz \cite[Thm.~1.3]{BlochEsnaultKerz2012}. Analogues in characteristic zero have also been established by them \cite{BlochEsnaultKerz2013} and the author \cite{Morrow_Deformational_Hodge}. 

We finish this introduction with a brief discussion of the proofs of Theorems \ref{theorem_3} and~\ref{theorem_4}. The new input which makes these results possible is recent work of the author joint with B.~Dundas \cite{Morrow_Dundas}, which establishes that topological cyclic homology is continuous under very mild hypotheses\footnote{(Added October 2014) I have a new proof of Theorem \ref{theorem_3}, hence also of Theorem \ref{theorem_2}, avoiding topological cyclic homology and higher algebraic $K$-theory; it will be included in forthcoming work on pro crystalline Chern characters and de Rham--Witt sheaves of infinitesimal thickenings.}. In particular, in the framework of Theorem \ref{theorem_4}, our results yield a weak equivalence \begin{equation}\TC(X;p)\stackrel{\sim}{\To}\op{holim}_r\TC(Y_r;p)\label{equ_we}\end{equation} between the $p$-typical topological cyclic homologies of $X$ and the limit of those of all the thickenings of the special fibre. Moreover, R.~McCarthy's theorem and the trace map describe the obstruction to infinitesimally lifting elements of $K$-theory in terms of $\TC(Y_r;p)$ and $\TC(Y;p)$, while results of T.~Geisser and L.~Hesselholt describe $\TC(X;p)$ and $\TC(Y;p)$ in terms of logarithmic de Rham--Witt cohomology since $X$ and $Y$ are regular. Combining these two descriptions with the weak equivalence (\ref{equ_we}) yields a preliminary form of Theorem \ref{theorem_4} phrased in terms of Gros' logarithmic crystalline Chern character; see Proposition \ref{proposition_preliminary_main}.

To deduce Theorem \ref{theorem_4} we then analyse the behaviour of the Frobenius on the crystalline cohomology of the $k[[t_1,\dots,t_m]]$-scheme $X$; see Proposition \ref{proposition_key}. Theorem \ref{theorem_3} then follows from Theorem \ref{theorem_4} using Grothendieck's algebrisation isomorphism $\Pic X\isoto\projlim_r\Pic Y_r$.

\subsection*{Acknowledgments}
I am very grateful to P.~Scholze for many discussions during the preparation of this article; in particular, his suggestion that Theorem \ref{theorem_3} be reformulated in terms of crystals and flatness has played a key role in the development of the article. It is also a pleasure to thank B.~Dundas for the collaboration which produced \cite{Morrow_Dundas}, the main results of which were proved with the present article in mind.

Moreover, I would like to thank D.~Caro, Y.~Nakkajima, A.~Shiho, and F.~Trihan for explaining aspects of arithmetic $\cal D$-modules, crystalline weights, log crystalline cohomology, and rigid cohomology to me, and S.~Bloch, H.~Esnault, and M.~Kerz for their interest in my work and their suggestions. The possibility of Corollary \ref{intro_corol_tate} was indicated to me by J.~Ayoub.

The existence of an error in the first version of the paper, involving line bundles with $\bb Q_p$-coefficients, was brought to my attention by F.~Charles.

\section{Conjecture \ref{VTC}}\label{section_general_comments}
Let $k$ be a perfect field of characteristic $p>0$, and write $W:=W(k)$, $K:=\Frac W(k)$. For any $k$-variety $X$, we denote by $H^n_\sub{crys}(X/W)$ and $H^n_\sub{crys}(X):=H^n_\sub{crys}(X/W)\otimes_WK$ its integral and rational crystalline cohomology groups. The crystalline cohomology $H^n_\sub{crys}(X)$ is naturally acted on by the absolute Frobenius $\phi:x\mapsto x^p$, whose $p^i$-eigenspaces we will denote by adding the traditional superscript $^{\phi=p^i}$. There is a theory of crystalline Chern classes, cycle classes, and the associated Chern character, constructed originally by P.~Berthelot and L.~Illusie \cite{BerthelotIllusie1970}, A.~Ogus \cite{Ogus1984}, and H.~Gillet and W.~Messing \cite{GilletMessing1987}:
\begin{align*}
&c_i(E)=c_i^\sub{crys}(E)\in\CH^i(X)_\bb Q\quad\mbox{($E$ a vector bundle on $X$),}\\
&cl_i=cl_i^\sub{crys}:\CH^i(X)_\bb Q\To H^{2i}_\sub{crys}(X),\\
&ch=ch^\sub{crys}:K_0(X)_\bb Q\To {\bigoplus}_{i\ge0}H^{2i}_\sub{crys}(X).
\end{align*}
The cycle classes and Chern character land in the $p^i$-eigenspace of $\phi$ acting on $H^{2i}_\sub{crys}(X)$.

For the reader's convenience, we restate the main conjecture from the Introduction;  the possible equivalence of (deform) and (crys) was raised by de Jong \cite{deJong2011} in a slightly different setting:

\setcounter{section}{0}
\begin{conjecture}[Crystalline Variational Tate Conjecture]\label{VTC}
Let $f:X\to S$ be a smooth, proper morphism of smooth $k$-varieties, $s\in S$ a closed point, and $z\in \CH^i(X_s)_\bb Q$. Let $c:=cl_i(z)\in H^{2i}_\sub{crys}(X_s)$. Then the following are equivalent:
\begin{description}
\item[(deform)] There exists $\tilde z\in \CH^i(X)_\bb Q$ such that $cl_i(\tilde z)|_{X_s}=c$.
\item[(crys)] $c$ lifts to $H^{2i}_\sub{crys}(X)$.
\item[(crys$\boldsymbol{-\phi}$)] $c$ lifts to $H^{2i}_\sub{crys}(X)^{\phi=p^i}$.
\item[(flat)] $c$ is flat, i.e., it lifts to $H^0_\sub{crys}(S,R^{2i}f_*\roi_{X/K})$.
\end{description}
\end{conjecture}

\setcounter{section}{1}\setcounter{theorem}{0}

More generally, we will discuss the conditions (crys), (crys-$\phi$), and (flat) for arbitrary cohomology classes $c\in H^{2i}_\sub{crys}(X_s)$.

\begin{remark}[The condition (flat)]\label{remark_flat}
The flatness condition, involving the global sections of Ogus' convergent $F$-isocrystal $R^{2i}f_*\roi_{X/K}$ \cite[\S3]{Ogus1984}, may require further explanation. The restriction map $H^{2i}_\sub{crys}(X)\to H^{2i}_\sub{crys}(X_s)$ may be factored as \[H^{2i}_\sub{crys}(X)\To H^0_\sub{crys}(S,R^{2i}f_*\roi_{X/K})\To H^{2i}_\sub{crys}(X_s),\] where the first arrow is an edge map in the Leray spectral sequence \[E^{ab}_2=H^a_\sub{crys}(S,R^bf_*\roi_{X/K})\implies H^{a+b}_\sub{crys}(X).\] We call an element $c\in H^{2i}_\sub{crys}(X_s)$ {\em flat} if and only if it lifts to $H^0_\sub{crys}(S,R^{2i}f_*\roi_{X/K})$; the lift, if it exists, is actually unique, assuming $S$ is connected \cite[Thm.~4.1]{Ogus1984}. In particular, the lift of a flat cycle class automatically lies in the eigenspace $H^0_\sub{crys}(S,R^{2i}f_*\roi_{X/K})^{\phi=p^i}$, so there is no need to introduce a (flat-$\phi$) condition. 
\end{remark}

\begin{remark}[Rigid cohomology 1]\label{remark_rigid_1}
In this remark we explain why we have chosen not to include Berthelot's rigid cohomology \cite{Berthelot1986} in the statement of Conjecture \ref{VTC}, even though it is {\em a priori} reasonable to consider the following conditions:
\begin{description}\itemsep1pt
\item[(rig)] $c$ lifts to $H^{2i}_\sub{rig}(X)$.
\item[(rig-$\boldsymbol{\phi}$)] $c$ lifts to $H^{2i}_\sub{rig}(X)^{\phi=p^i}$.
\end{description}
Indeed, by the obvious implications which we will mention in Remark \ref{remark_obvious}, the validity of Conjecture \ref{VTC} is unchanged by the addition of conditions (rig) and (rig-$\phi$). 

The rigid analogue of condition (flat) is more subtle, but also redundant, as we now explain. Let $j^*:F\mbox{-Isoc}^\dag(S/K)\to F\mbox{-Isoc}(S/K)$ denote the forgetful functor from overconvergent $F$-isocrystals on $S$ to convergent $F$-isocrystals on $S$. It is an open conjecture of Berthelot \cite[\S4.3]{Berthelot1986} that Ogus' convergent $F$-isocrystal $R^{2i}f_*\roi_{X/K}$ admits an overconvergent extension; that is, there exists $R^{2i}f_*\roi_{X/K}^\dag\in F\mbox{-Isoc}^\dag(S/K)$ such that $j^*(R^{2i}f_*\roi_{X/K}^\dag)=R^{2i}f_*\roi_{X/K}$. Moreover, the functor $j^*$ is now known to be fully faithful thanks to K.~Kedlaya \cite{Kedlaya2004}, and so $R^{2i}f_*\roi_{X/K}^\dag$ is unique if it exists.

Assuming for a moment the validity of Berthelot's conjecture, the following rigid analogue of (flat) could be considered:
\begin{description}
\item[(rig-flat)] $c$ lifts to $H^0_\sub{rig}(S,R^{2i}f_*\roi_{X/K}^\dag)$.
\end{description}
However, the conditions (flat) and (rig-flat) are in fact equivalent for any $c\in H^{2i}_\sub{crys}(X_s)^{\phi=p^i}$. Firstly, the implication $\Leftarrow$ is trivial. Secondly, assuming $c$ lifts to $\tilde c\in H^0_\sub{crys}(S,R^{2i}f_*\roi_{X/K})$, the lift $\tilde c$ automatically belongs to $H^0_\sub{crys}(S,R^{2i}f_*\roi_{X/K})^{\phi=p^i}$, as explained at the end of Remark \ref{remark_flat}. But, in the following diagram,
\[\xymatrix@R=2mm{
H^0_\sub{rig}(S,R^{2i}f_*\roi_{X/K}^\dag)^{\phi=p^i}\ar@{=}[d]\ar[r] & H^0_\sub{crys}(S,R^{2i}f_*\roi_{X/K})^{\phi=p^i}\ar@{=}[d]\\
\Hom_{F\mbox{-Isoc}^\dag(S/K)}(\roi_{S/K}^\dag(i),R^{2i}f_*\roi_{X/K}^\dag)\ar[r] & \Hom_{F\mbox{-Isoc}(S/K)}(\roi_{S/K}(i),R^{2i}f_*\roi_{X/K}),
}\]
where $(i)$ denotes Tate twists and all other notation should be clear, the bottom horizontal arrow is an isomorphism by Kedlaya's aforementioned fully faithfulness result. Therefore $\tilde c$ lifts uniquely to $H^0_\sub{rig}(S,R^{2i}f_*\roi_{X/K}^\dag)^{\phi=p^i}$, proving the implication $\Rightarrow$.

In conclusion, since (rig-flat) is only well-defined conditionally on Berthelot's conjecture, and since it is then equivalent to (flat), we will not consider it further; the exception is Remark \ref{remark_rigid_2}, where we make further comments on rigid cohomology.
\end{remark}

\begin{remark}[Obvious implications]\label{remark_obvious}
Since rigid cohomology maps to crystalline cohomology (e.g., via the inclusion $W^\dag\Omega_X^\blob\subseteq W\Omega_X^\blob$ of Davis--Langer--Zink's overconvergent de Rham--Witt complex into the usual de Rham--Witt complex \cite{LangerZinkDavis2011}), we have the following automatic implications for any $c\in H^{2i}_\sub{crys}(X_s)$:
\[\xymatrix@=5mm{
{\bf(rig\text-\boldsymbol{\phi})}\ar@{=>}[r]\ar@{=>}[d] & {\bf(crys\text-\boldsymbol{\phi})}\ar@{=>}[d]& \\
{\bf(rig)}\ar@{=>}[r] & {\bf(crys)}\ar@{=>}[r]&{\bf(flat)}
}\]
Furthermore, if $c=cl_i(z)$ for some $z\in \CH^i(X_s)_\bb Q$, then these five conditions on $c$ are all consequences of the condition (deform), since the cycle class map $cl_i:\CH^i(X)_\bb Q\to H^{2i}_\sub{crys}(X)$ factors through $H^{2i}_\sub{rig}(X)^{\phi=p^i}$ by \cite{Petrequin2003}.
\end{remark}

Our main result is the proof of Conjecture \ref{VTC} for divisors (henceforth identified with line bundles), assuming $f$ is projective. We finish this section by proving this, assuming the main results of later sections, and then presenting a variant for hypersurfaces:

\begin{theorem}\label{theorem_VTC_for_divisors}
Let $f:X\to S$ be a smooth, projective morphism of smooth $k$-varieties, $s\in S$ a closed point, and $L\in\Pic(X_s)_\bb Q$. Let $c:=c_1(L)\in H^{2}_\sub{crys}(X_s)$. Then the following are equivalent:
\begin{description}
\item[(deform)] There exists $\tilde L\in \Pic(X)_\bb Q$ such that $c_1(\tilde L)|_{X_s}=c$.
\item[(crys)] $c$ lifts to $H^{2}_\sub{crys}(X)$.
\item[(crys$\boldsymbol{-\phi}$)] $c$ lifts to $H^{2}_\sub{crys}(X)^{\phi=p}$.
\item[(flat)] $c$ is flat, i.e., it lifts to $H^0_\sub{crys}(S,R^{2}f_*\roi_{X/K})$.
\end{description}
\end{theorem}
\begin{proof}[Proof (assuming Corol.~\ref{corollary_line_bundles} and Thm.~\ref{theorem_Partie_Fixe})]
We may assume $S$ is connected. By Theorem \ref{theorem_Partie_Fixe} below and the obvious implication (deform)$\Rightarrow$(crys), it is enough to prove the implication (crys-$\phi$)$\Rightarrow$(deform). So assume that $c$ lifts to $\tilde c\in H^{2}_\sub{crys}(X)^{\phi=p}$. Let $A:=\roi_{S,s}^\sub{sh}$ be the strict Henselisation of $\roi_{S,s}$, and let $L^\sub{alg}$ denote the pullback of $L$ to $X\times_Sk(s)^\sub{alg}$.

Applying Corollary \ref{corollary_line_bundles} to $\hat X:=X\times_S\hat A$ we see that there exists $L_1\in\Pic(\hat X)[\tfrac1p]$ such that of $L_1|_{X\times_Sk(s)^\sub{alg}}=L^\sub{alg}$. The rest of the proof consists of descending $L_1$ from $\hat X$ to $X$.

By N\'eron--Popescu desingularisation \cite{Popescu1985, Popescu1986}, we may write $\hat A$ as a filtered colimit of smooth, local $A$-algebras. Therefore there exist a smooth, local $A$-algebra $A'$, a morphism of $A$-algebras $A'\to\hat A$, and a line bundle $L_2\in \Pic(X\times_S A')[\tfrac1p]$ such that $L_2$ pulls back to $L_1$ via the morphism $\hat X\to X\times_S A'$. As usual, the composition $A\to A'\to\hat A$ induces an isomorphism of residue fields, so that $A\to A'$ has a section at the level of residue fields; but $A$ is Henselian and $A\to A'$ is smooth, so this lifts to a section $\sigma:A'\to A$. Then, by construction, the restriction of $L_3:=\sigma^*L_2\in\Pic(X\times_SA)[\tfrac1p]$ to $X\times_Sk(s)^\sub{alg}$ is $L^\sub{alg}$.

But $A$ is the filtered colimit of the connected \'etale neighbourhoods of $\Spec k(s)^\sub{alg}\to S$; so there exists an \'etale morphism $U\to S$ (with $U$ connected), a closed point $s'\in U$ sitting over $s$, and a line bundle $L_4\in\Pic(X\times_SU)[\tfrac1p]$ such that the restriction of $L_4$ to $X\times_Sk(s')$ coincides with the pullback of $L$ to $X\times_Sk(s')$.

Let $S'$ be the normalisation of $S$ inside the function field of $U$, and let $L_5\in\Pic(X\times_SS')[\tfrac1p]$ be any extension of $L_4$, which exists by normality of $X\times_SS'$. By de Jong's theory of alterations \cite{deJong1996}, there exists a generically \'etale alteration $\pi'':S''\to S'$ with $S''$ connected and smooth over $k$. Let $L_6:=\pi''^*L_5\in\Pic(X\times_SS'')[\tfrac1p]$, and also let $s''\in S''$ be any closed point sitting over $s'$. To summarise, we have a commutative diagram

\[\xymatrix{
X \ar[d]^f & X' \ar[l]_{\pi'_X}\ar[d]^{f'}& X''\ar[l]_{\pi''_X}\ar[d]^{f''}\\
S & S' \ar[l]_{\pi'} & S'' \ar[l]_{\pi''} \\
\Spec k(s)\ar@{^(->}[u]&\Spec k(s')\ar@{^(->}[u]\ar[l]&\Spec k(s'')\ar[l]\ar@{^(->}[u]&
}\]
where:
\begin{itemize}\itemsep0pt
\item $X':=X\times_SS'$ and $X'':=X\times_SS''$;
\item $\pi:=\pi'\circ\pi''$ is a generically \'etale alteration;
\item the restriction of $L_6$ to $X''\times_{S''}k(s'')=X\times_Sk(s'')$ coincides with the pullback to $L$ to $X\times_Sk(s'')$.
\end{itemize}
By studying crystalline Chern classes, we can now complete the proof. It will be convenient to denote by $e:H^2_\sub{crys}(X)\to H^0_\sub{crys}(S,R^2f_*\roi_{X/K})$ the edge map in the Leray spectral sequence (abusing notation, we also use the notation $e$ for the families $X'$, $X''$, etc.); in particular, set $\res c:=e(\tilde c)$.

Let $V\subseteq X$ be a nonempty open subscheme such that $\pi^{-1}(V)\to V$ is finite \'etale (this exists since $\pi$ is proper and generically \'etale), and let $L_7\in\Pic(X\times_SV)[\tfrac1p]$ be the pushforward of $L_6|_{X\times_S\pi^{-1}(V)}$. We claim that $e(c_1(L_7))=n\,\res c|_V$, where $n$ is the generic degree of $\pi$.

First we note that $e(c_1(L_6))=\pi^*(\res c)$ in $H^0_\sub{crys}(S'',R^2f''_*\roi_{X''/K})$: the two classes agree at $s''$ by construction of $L_6$, and the specialisation map $H^0_\sub{crys}(S'',R^2f''_*\roi_{X''/K})\to H^2_\sub{crys}(X''\times_{S''}k(s''))$ is injective by \cite[Thm.~4.1]{Ogus1984}. Restricting to $\pi^{-1}(V)$, we therefore obtain $e(c_1(L_6|_{X\times_S\pi^{-1}(V)}))=\pi^*(\res c|_V)$. Pushing forwards along the finite \'etale morphism $\pi^{-1}(V)\to V$ proves the claim.

Finally, let $\tilde L\in\Pic(X)_\bb Q$ be any extension of $L_7^{1/n}$. Then $e(c_1(\tilde L))$ agrees with $\res c$ on $V$, hence agrees with $\res c$ everywhere. In particular, by specialising to $s$ we obtain that $c_1(\tilde L)|_{X_s}=c$; this completes the proof.
\end{proof}

The following consequence of the main result concerns hyperplane sections:

\begin{corollary}\label{corollary_hypersurfaces}
Let $X$ be a smooth, projective $k$-variety, $Y\into X$ a smooth ample divisor, and $L\in\Pic(Y)_\bb Q$. Let $c:=c_1(L)\in H^2_\sub{crys}(Y)$. Then the following are equivalent:
\begin{enumerate}
\item There exists $\tilde L\in\Pic(X)_\bb Q$ such that $c_1(\tilde L)|_Y= c$.
\item $c$ lifts to $H^2_\sub{crys}(X)$.
\item $c$ lifts to $H^2_\sub{crys}(X)^{\phi=p}$.
\end{enumerate}
\end{corollary}
\begin{proof}
We begin by fitting $Y$ into a pencil of hyperplane sections in the usual way. Choosing a line in the linear system $|\roi(Y)|$ yields a pencil of hyperplane sections $\{X_t:t\in\bb P^1_k\}$, with base locus denoted by $B\into X$, and an associated fibration; that is, there is a diagram \[X\stackrel{\pi}{\longleftarrow} X'\stackrel{f}{\To}\bb P^1_k,\] where:
\begin{itemize}\itemsep0pt
\item $\pi$ is the blow-up of $X$ along the smooth subvariety $B$.
\item $f$ is projective and flat, and is smooth over a non-empty open $V\subseteq\bb P_k^1$ which contains $0$.
\item For each closed point $t\in \bb P_k^1$, the fibre $X'_t$ is isomorphic via $\pi$ to the hyperplane section $X_t$ of $X$; in particular, $X_0'\cong X_0=Y$.
\end{itemize}
We may now properly begin the proof; the only non-trivial implication is (ii)$\Rightarrow$(i), so assume that $c$ lifts to $H^2_\sub{crys}(X)$. Then $c$ certainly lifts to $H^2_\sub{crys}(f^{-1}(V))$, so Theorem \ref{theorem_VTC_for_divisors}, applied to the smooth, projective morphism $f^{-1}(V)\to V$, implies that there exists $L'\in\Pic(f^{-1}(V))_\bb Q$ such that $c^1(L')|_{X'_0}=c$. Then $L'$ may be spread out to some $L''\in\Pic(X')_\bb Q$, which evidently still satisfies $c_1(L'')|_{X'_0}=c$.

Let $E:=\pi^{-1}(B)$ denote the exception divisor of the blow-up $\pi$. By the standard formula for the Picard group of a blow-up along a regularly embedded subvariety, we may write $L''=\pi^*(L''')\otimes\roi(E)^a$ for some $L'''\in\Pic(X)_\bb Q$ and $a\in\bb Q$. But the line bundles $\roi(E)$ and $\pi^*(\roi(Y))$ have the same restriction to $Y=X_0\cong X_0'$, namely $\roi(B)$. Hence $\pi^*(L'''\otimes\roi(Y)^a)$ has the same restriction to $Y$ as $L''$; so, setting $\tilde L:=L'''\otimes\roi(Y)^a\in\Pic(X)_\bb Q$, we see that $c_1(\tilde L)|_Y=c_1(L'')|_{X'_0}=c$, as desired.
\end{proof}

\begin{remark}\label{remark_hyerplane_theorem}
The only interesting case of Corollary \ref{corollary_hypersurfaces} is when $\dim X=3$. Indeed, if $\dim X>3$ then $H^2_\sub{crys}(X)\isoto H^2_\sub{crys}(X_t)$ by Weak Lefschetz for crystalline cohomology \cite[\S3.8]{Illusie1975} and $\Pic(X)_\bb Q\isoto\Pic(Y)_\bb Q$ by Grothendieck--Lefschetz \cite[Exp.~XI]{SGA_II}. On the other hand, if $\dim X=2$ then $H^2(Y)$ is one-dimensional, spanned by $c_1(\roi(Y))|_Y$, and the line bundle $\tilde L$ may always be taken to be $\roi(Y)^{\deg L}$.

When $\dim X=3$, the restriction maps $H^2_\sub{crys}(X)\to H^2_\sub{crys}(Y)$ and $\Pic(X)_\bb Q\to\Pic(Y)_\bb Q$ are injective by Weak Lefschetz and Grothendieck--Lefschetz respectively, so Corollary \ref{corollary_hypersurfaces} may then be more simply stated as follows: the line bundle $L\in\Pic(Y)_\bb Q$ lies in $\Pic(X)_\bb Q$ if and only if its Chern class $c_1(L)\in H^2_\sub{crys}(Y)$ lies in $H^2_\sub{crys}(X)$.
\end{remark}

\section{Crystalline Th\'eor\`eme de la Partie Fixe}\label{section_partie_fixe}
The aim of this section is to prove the following equivalences between the conditions appearing in Conjecture \ref{VTC} ($k$ continues to be a perfect field of characteristic $p$):

\begin{theorem}\label{theorem_Partie_Fixe}
Let $f:X\to S$ be a smooth, projective morphism of smooth $k$-varieties, $s\in S$ a closed point, and $c\in H^{2i}_\sub{crys}(X_s)^{\phi=p^i}$. Then the following are equivalent:
\begin{description}
\item[(crys)] $c$ lifts to $H^{2i}_\sub{crys}(X)$.
\item[(crys-$\boldsymbol{\phi}$)] $c$ lifts to $H^{2i}_\sub{crys}(X)^{\phi=p^i}$.
\item[(flat)] $c$ is flat, i.e., it lifts to $H^0_\sub{crys}(S,R^{2i}f_*\roi_{X/K})$.
\end{description}
\end{theorem}

\begin{remark}
Theorem \ref{theorem_Partie_Fixe} and its proof via Theorem \ref{theorem_crystalline_partie_fixe} resemble P.~Deligne's Th\'eor\`eme de la Partie Fixe \cite[\S4.1]{Deligne1971} for de Rham cohomology in characteristic zero, though the name in our case is a misnomer as we do not consider the action of the fundamental group $\pi_1(S,s)$.

Deligne extends his result in characteristic zero to smooth, {\em proper} morphisms using resolution of singularities. Unfortunately, the standard arguments with alterations appear to be insufficient for us to do the same in characteristic $p$, and so we are forced to restrict to projective morphisms in some of our main results.
\end{remark}

We begin with two preliminary results, Lemmas \ref{lemma_SS} and \ref{lemma_SS_2}, on spectral sequences in an arbitrary abelian category. In these lemmas all spectral sequences are implicitly assumed to start on the $E_1$-page for simplicity. We say that the $n^\sub{th}$ column of a spectral sequence $E_*^{ab}$ is {\em stable} if and only if $E_1^{n b}=E_\infty^{n b}$ for all $b\in\bb Z$; in other words, all differentials into and out of the $n^\sub{th}$ column are zero. Obvious modification of the terminology, such as {\em stable in columns $> n$}, or $<n$, will be used in this section.

\begin{lemma}\label{lemma_SS}
Suppose that $E_*^{ab}$ and $F_*^{ab}$ are spectral sequences, that $n\in\bb Z$, and that the following conditions hold:
\begin{enumerate}\itemsep0pt
\item The spectral sequence $E_*^{ab}$ is stable in columns $>n$.
\item The spectral sequence $F_*^{ab}$ is stable in columns $<n$.
\item There is a map of spectral sequences $f:E_*^{ab}\to F_*^{ab}$ which is an isomorphism on the $n^\sub{th}$ columns of the first pages, i.e., $f:E_1^{n b}\isoto F_1^{n b}$ for all $b\in\bb Z$.
\end{enumerate}
Then the $n^\sub{th}$ columns of both spectral sequences are stable, and isomorphic via $f$, i.e.,
\[\xymatrix@R=4mm{
E_1^{n b} \ar[r]^{f\,\cong}\ar@{=}[d] & F_1^{n b}\ar@{=}[d]\\
E_\infty^{n b} \ar[r]_{f\,\cong} & F_\infty^{n b}
}\]
\end{lemma}
\begin{proof}
According to assumptions (i) and (ii), all differentials with domain (resp.~codomain) in the $n^\sub{th}$ column of any page of the $E$-spectral sequence (resp.~$F$-spectral sequence) are zero. So it remains to check that all differentials with codomain (resp.~domain) in the $n^\sub{th}$ column of any page of the $E$-spectral sequence (resp.~$F$-spectral sequence) are zero. This is an easy induction, using assumption (iii), on the page number of the spectral sequence.
\end{proof}

The following technique to check the degeneration of a family of spectral sequences is inspired by \cite[Thm.~1.5]{Deligne1968}:

\begin{lemma}\label{lemma_SS_2}
Let $d\ge0$ and let \[\cdots\xto{u}E^{ab}_*(-4)\xto{u}E^{ab}_*(-2)\xto{u}E^{ab}_*(0)\xto{u}E^{ab}_*(2)\xto{u}E^{ab}_*(4)\xto{u}\cdots\] be a sequence of right half plane spectral sequences. Make the following assumptions:
\begin{enumerate}\itemsep0pt
\item For every $n\in 2\bb Z$, the spectral sequence $E^{ab}_*(n)$ vanishes in columns $>n$.
\item For every $n\in\bb Z$ and $i\ge0$ such that $i\equiv n$ mod $2$, the map of spectral sequences $u^i:E^{ab}_*(n-i)\to E^{ab}_*(n+i)$ is an isomorphism on the $n-d^\sub{th}$ columns of the first pages.
\end{enumerate}
Then, for every $n\in2\bb Z$, the spectral sequence $E_*^{ab}(n)$ degenerates on the $E_1$-page.
\end{lemma}
\begin{proof}
For any $n\in2\bb Z$ and any integer $c\le d+n$, assumption (ii) implies that \begin{equation}E_*^{ab}(2c-n)\xto{u^{d+n-c}}E_*^{ab}(2d+n)\label{equ_SS1}\end{equation} is an isomorphism on the $c^\sub{th}$ columns of the first pages. In particular, if $c<n$, so that the $c^\sub{th}$ column of the left spectral sequence in (\ref{equ_SS1}) vanishes, then the $c^\sub{th}$ column of the right spectral sequence also vanishes. That is, $E_*^{ab}(2d+n)$ vanishes in columns $<n$; or, reindexing, $E_*^{ab}(n)$ vanishes in columns $<n-2d$ for every $n\in2\bb Z$.

For each $i=0,\dots,d+1$, we now make the following claim: {\em for every $n\in2\bb Z$, the spectral sequence $E_*^{ab}(n)$ is stable in columns $>n-i$ and in columns $<n+i-2d$}.

The claim is true when $i=0$, thanks to assumption (i) and our vanishing observation above. Proceeding by induction, assume that the claim is true for some $i\in\{0,\dots,d\}$. Then, for any $n\in2\bb Z$, the map of spectral sequences \begin{equation}E_*^{ab}(n)\xto{u^{d-i}}E_*^{ab}(n+2d-2i)\label{equ_SS2}\end{equation} is an isomorphism on the $n-i$ columns of the $E_1$-pages, by assumption (ii); moreover, by the inductive hypothesis, the left spectral sequence in (\ref{equ_SS2}) is stable in columns $>n-i$ and the right spectral sequence in columns $<n-i$. By Lemma \ref{lemma_SS}, both spectral sequences are therefore stable in column $n-i$, proving the inductive claim for $i+1$.

But this completes the proof, for the inductive claim at $i=d+1$ asserts that the spectral sequence $E^{ab}_*(n)$ is stable in columns $>n-d-1$ and $<n-d+1$, hence degenerates on the first page.
\end{proof}

\begin{remark}[Hard Lefschetz for crystalline cohomology]\label{remark_Lef}
Before proving the main theorems of the section we make a remark on the Hard Lefschetz theorem for crystalline cohomology. Let $X$ be a smooth, projective, connected, $d$-dimensional variety over $k$, and let $L$ be an ample line bundle on $X$. Let $u:=c_1(L)\in H^2_\sub{crys}(X)$, and also denote by $u$ the induced cup product map $u\cup-:H^*_\sub{crys}(X)\to H^{*+2}_\sub{crys}(X)$. We understand the Hard Lefschetz theorem as the assertion that \begin{equation}u^i: H^{d-i}_\sub{crys}(X)\To H^{d+i}_\sub{crys}(X)\label{equ_HL}\end{equation} is an isomorphism of $K$-vector spaces for $i=0,\dots,d$.

Assuming in addition that $L=\roi(D)$ for some smooth hyperplane section $D$ of $X$, that $X$ is geometrically connected over $k$, and that $k$ is finite, isomorphism (\ref{equ_HL}) follows from the $\ell$-adic case, as explained in \cite{KatzMessing1974}. We now explain how to reduce the more general assertion above to this special case. Firstly, we may assume that $X$ is geometrically connected over $k$ by replacing $k$ by $H^0(X,\roi_X)$; then we may assume $k$ is finite by a standard spreading out argument \cite[\S3.8]{Illusie1975}; thirdly we may assume $L$ is very ample by replacing $L$ by $L^m$, as this merely replaces $u$ by $mu$; and finally we may assume $L=\roi(D)$ for some smooth hyperplane section $D$ of $X$, by passing to a finite extension of $k$ after which $L=i^*\roi(1)$ for some closed embedding $i:X\into\bb P^N_k$ such that there exists a hyperplane in $\bb P^N_k$ having smooth intersection with $X$.
\end{remark}

The following is our analogue in crystalline cohomology of Deligne's Th\'eor\`eme de la Partie Fixe \cite[\S4.1]{Deligne1971}, in which $\phi$ denotes, as everywhere, the absolute Frobenius:

\begin{theorem}\label{theorem_crystalline_partie_fixe}
Let $f:X\to S$ be a smooth, projective morphism of smooth $k$-varieties. Then:
\begin{enumerate}
\item The Leray spectral sequence $E^{ab}_2=H^a_\sub{crys}(S,R^bf_*\roi_{X/K})\Rightarrow H^{a+b}_\sub{crys}(X)$ degenerates on the $E_2$-page.
\item For any $r\ge0$, $s\in\bb Z$, the Leray spectral sequence in (i) contains a sub spectral sequence $H^a_\sub{crys}(S,R^bf_*\roi_{X/K})^{\phi^r=p^s}\Rightarrow H^{a+b}_\sub{crys}(X)^{\phi^r=p^s}$ which also degenerates on the $E_2$-page. (The superscripts denote the $p^s$-eigenspaces of $\phi^r$.)
\item For any $n,r\ge0$, $s\in\bb Z$, the canonical map \[H^n_\sub{crys}(X)^{\phi^r=p^s}\To H^0_\sub{crys}(S,R^nf_*\roi_{X/K})^{\phi^r=p^s}\] is surjective. 
\end{enumerate}
\end{theorem}
\begin{proof}
We may assume that $S$ and $X$ are connected, and we let $d$ denote the relative dimension of $f$. Let $L$ be a line bundle on $X$ which is relatively ample with respect to $f$, and let $u:=c_1(L)\in H^2_\sub{crys}(X)$. Also denote by $u$ the induced cup product morphism of convergent isocrystals $u:R^if_*\roi_{X/K}\to R^{i+2}f_*\roi_{X/K}$, and note that \begin{equation}u^i:R^{d-i}f_*\roi_{X/K}\To R^{d+i}f_*\roi_{X/K}\label{equ_HL2}\end{equation} is an isomorphism for $i\ge0$; indeed, it is enough to check this isomorphism after restricting to each closed point of $S$ \cite[Lem.~3.17]{Ogus1984}, where it is exactly (using the identification in \cite[Rem.~3.7.1]{Ogus1984}) the Hard Lefschetz theorem for crystalline cohomology which was discussed in Remark \ref{remark_Lef}. Note that isomorphism (\ref{equ_HL2}) is valid even if $i>d$, the left side being zero by convention and the right side by relative dimension considerations; this is helpful for indexing.

Claim (i) now follows by applying Deligne's axiomatic approach to the degeneration of Leray spectral sequences via Hard Lefschetz \cite[\S1]{Deligne1968} to $Rf_*\roi_{X/K}$, which lives in the derived category of Ogus' convergent topos $(S/W)_\sub{conv}^\sim$ \cite{Ogus1990}.

(ii): By (i), each group $H^n:=H^n_\sub{crys}(X)$ (which, for simplicity of indexing, is defined to be zero if $n<0$) is equipped with a natural descending filtration \[H^n=\cdots =F_{-1}H^n=F_0H^n\supseteq \cdots\supseteq F_nH^n\supseteq F_{n+1}H^n=F_{n+2}H^n=\cdots=0\] having graded pieces $g_a^n=F_aH^n/F_{a+1}H^n\cong H^a_\sub{crys}(S,R^{n-a}f_*\roi_{X/K})$ for $a\in\bb Z$. Note that this filtration is respected by $\phi$. The assertion to be proved is that the induced filtration on the subgroup $(H^n)^{\phi^r=p^s}$ has graded pieces $(g_a^n)^{\phi^r=p^s}$.

In other words, fixing $r\ge0$ but allowing $s\in\bb Z$ to vary, we must show that the Ker--Coker spectral sequence for the morphism $\phi^r-p^s:H^n\to H^n$ of filtered groups, namely \[E_1^{ab}(n,s)=\begin{cases}\ker(g_a^n\xto{\phi^r-p^s}g_a^n) & a+b=0 \\ \op{Coker}(g_a^n\xto{\phi^r-p^s}g_a^n) & a+b=1 \\ 0 & \mbox{else}\end{cases}\implies \begin{cases}\ker (H^n\xto{\phi^r-p^s}H^n) & a+b=0 \\ \op{Coker} (H^n\xto{\phi^r-p^s}H^n) & a+b=1 \\ 0&\mbox{else}\end{cases}\] degenerates on the first page.

Since $\phi u=pu\phi$, cupping with $u$ defines a morphism of spectral sequences $u:E^{ab}_*(n,s)\to E^{ab}_*(n+2,s+1)$, and we will check that this family of spectral sequences satisfies the hypotheses of Lemma \ref{lemma_SS_2} (to be precise, we apply Lemma \ref{lemma_SS_2} to $E^{ab}_*(n):=\bigoplus_{s\in\bb Z}E^{ab}_*(n,s)$, where $n\in 2\bb Z$ and with $u$ defined degree-wise). Firstly, $E^{ab}_*(n,s)$ vanishes in columns $>n$ since the filtration on $H^n$ has length $n$. Secondly, the morphism \[u^i:E^{ab}_*(n-i,s)\to E^{ab}_*(n+i,s+i)\] is an isomorphism on the $n-d^\sub{\,th}$ column of the $E_1$-pages for any $n\in\bb Z$ by (\ref{equ_HL2}). Hence, by Lemma \ref{lemma_SS_2}, the spectral sequence $E^{ab}_*(n,s)$ degenerates on the first page for every $n\in 2\bb Z$ and $s\in\bb Z$; a minor reindexing treats the case that $n$ is odd, completing the proof of (ii).

Claim (iii) is immediate from (ii), as the canonical map is the edge map.
\end{proof}

\begin{remark}
More generally, Theorem \ref{theorem_crystalline_partie_fixe}(ii)\&(iii) remain true if $\phi^r-p^s$ is replaced by any $K$-linear combination of integral powers of $\phi$.
\end{remark}

\begin{remark}
Deligne's result used in the proof of Theorem \ref{theorem_crystalline_partie_fixe}(i) states not only that the Leray spectral sequence degenerates, but even that there is a non-canonical isomorphism $Rf_*\roi_{X/K}\cong\bigoplus_iR^if_*\roi_{X/K}[-i]$ in $D^b((S/W)_\sub{conv}^\sim)$. The content of Lemmas \ref{lemma_SS} and \ref{lemma_SS_2} is essentially that this isomorphism may be chosen to be compatible with the action of $\phi$.
\end{remark}

Theorem \ref{theorem_crystalline_partie_fixe} is evidently sufficient to prove our desired equivalences:

\begin{proof}[Proof of Theorem \ref{theorem_Partie_Fixe}]
We may assume $S$ is connected. In light of the obvious implications of Remark \ref{remark_obvious}, it is enough to prove that (flat)$\Rightarrow$(crys-$\phi$); so assume that $c$ lifts to some $\tilde c\in H^0_\sub{crys}(S,R^{2i}f_*\roi_{X/K})$. As discussed in Remark \ref{remark_flat}, the canonical map $H^0_\sub{crys}(S,R^{2i}f_*\roi_{X/K})\to H^{2i}_\sub{crys}(X_s)$ is injective by \cite[Thm.~4.1]{Ogus1984}, and hence $\tilde c\in H^0_\sub{crys}(S,R^{2i}f_*\roi_{X/K})^{\phi=p^i}$. Theorem \ref{theorem_crystalline_partie_fixe}(iii) implies that $\tilde c$ lifts to $H^{2i}_\sub{crys}(X)^{\phi=p^i}$, completing the proof.
\end{proof}

\begin{remark}[Rigid cohomology 2]\label{remark_rigid_2}
We finish this section by making a further remark on rigid cohomology, continuing Remark \ref{remark_rigid_1}. We assume throughout this remark that $f:X\to S$ is a smooth, proper morphism, where $S$ is a smooth, affine curve over $k=k^\sub{alg}$, and that $f$ admits a semi-stable compactification $\res f:\res X\to\res S$, where $\res S$ and $\res X$ are smooth compactifications of $S$ and $X$. We will sketch a proof of the following strengthening (which we do not need) of Theorems \ref{theorem_Partie_Fixe} and \ref{theorem_crystalline_partie_fixe} in this special case:
\begin{quote}
The composition \[H^{2i}_\sub{rig}(X)^{\phi=p^i}\To H^{2i}_\sub{crys}(X)^{\phi=p^i}\To H^0_\sub{crys}(S,R^{2i}f_*\roi_{X/K})^{\phi=p^i}\] is surjective.
\end{quote}

Since $S$ is an affine curve, Berthelot's conjecture discussed in Remark \ref{remark_rigid_1} is known to have an affirmative answer by S.~Matsuda and F.~Trihan \cite{TrihanMatsuda2004}: a unique overconvergent extension $R^if_*\roi_{X/K}^\dag$ of $R^if_*\roi_{X/K}$ exists for each $i\ge0$. It is then not unreasonable to claim that there ``obviously'' exists a Leray spectral sequence in rigid cohomology, namely \begin{equation}E_2^{ab}=H^a_\sub{rig}(S,R^bf_*\roi_{X/K}^\dag)\implies H^{a+b}_\sub{rig}(X),\label{rigid_Leray}\end{equation} but this turns out to be highly non-trivial. Assuming for a moment that this spectral sequence exists, it must degenerate for dimension reasons and thus yield short exact sequences \[0\To H^1_\sub{rig}(S,R^{n-1}f_*\roi_{X/K}^\dag)\To H^n_\sub{rig}(X)\To H^0_\sub{rig}(S,R^nf_*\roi_{X/K}^\dag)\To 0.\] By finiteness of rigid cohomology these groups are all finite dimensional $K$-vector spaces, and so by Dieudonn\'e--Manin the sequence remains exact after restricting to Frobenius eigenspaces. In particular, the edge maps \[H^{2i}_\sub{rig}(X)^{\phi=p^i}\To H^0_\sub{rig}(S,R^{2i}f_*\roi_{X/K}^\dag)^{\phi=p^i}=H^0_\sub{crys}(S,R^{2i}f_*\roi_{X/K})^{\phi=p^i}\] (the final equality was explained in Remark \ref{remark_rigid_1}) are surjective, as desired.

It remains to show that the rigid Leray spectral sequence (\ref{rigid_Leray}) really exists; the proof uses results from log crystalline cohomology due to A.~Shiho. Let $M=\res S\setminus S$ and $D=\res X\setminus X$, let $(\res S,M)$ and $(\res X,D)$ denote the associated log schemes, and view $\res f:(\res X,D)\to(\res S,M)$ as a log smooth morphism of Cartier type. By the general theory of sites there is an associated Leray spectral sequence in log-convergent cohomology \begin{equation}E_2^{ab}=H^a((\res S,M)_\sub{logconv},R^bf_*\cal K)\implies H^{a+b}((\res X,D)_\sub{logconv},\cal K),\label{log_Leray}\end{equation} where $\cal K$ is the structure sheaf on $(\res X,D)_\sub{logconv}$, the log-convergent site of $(\res X,D)$. However, \cite[Corol.~2.3.9 \& Thm.~2.4.4]{Shiho2002} state that $H^{a+b}((\res X,D)_\sub{logconv},\cal K)\cong H^{a+b}_\sub{rig}(X)$. It can moreover be shown, using Shiho's results on relative log-convergent cohomology \cite{Shiho2007I,Shiho2007II,Shiho2008}, that the terms on the $E_2$-page of (\ref{log_Leray}) are naturally isomorphic to $H^a_\sub{rig}(S,R^bf_*\roi_{X/K}^\dag)$, thereby completing the proof of the existence of (\ref{rigid_Leray}). This completes the sketch of the proof of the strengthening of Theorems \ref{theorem_Partie_Fixe} and \ref{theorem_crystalline_partie_fixe} in this special case.
\end{remark}

\section{Local and infinitesimal forms of Conjecture \ref{VTC}}
The aim of this section is to prove Theorems \ref{theorem_3} and \ref{theorem_4} from the Introduction, which are local and infinitesimal analogues of Conjecture \ref{VTC}. As always, $k$ is a perfect field of characteristic $p$. We will study a smooth, proper scheme $X$ over $\Spec A$, where $A:=k[[t_1,\dots,t_m]]$, and we adopt the following notation: the special fibre and its infinitesimal thickenings inside $X$ are always denoted by \[Y:=X\times_Ak, \qquad Y_r:=X\times_AA/\pid{t_1,\dots,t_m}^r.\]
Some of the deformation results of this section only require $p$-torsion to be neglected, so we write $M[\tfrac1p]:=M\otimes_\bb Z\bb Z[\tfrac1p]$ for any abelian group $M$.

\subsection{Some remarks on crystalline cohomology}\label{subsection_remarks}
The proofs of Theorems \ref{theorem_3} and \ref{theorem_4} require us to work with crystalline cohomology and the de Rham--Witt sheaves $W_r\Omega^n_X$ of S.~Bloch, P.~Deligne, and L.~Illusie, in greater generality than can be found in Illusie's treatise \cite{Illusie1979}; some of the generalisations we need can be readily extracted from his proofs and have already been presented in \cite{Shiho2007}, while others follow via a filtered colimit argument to reduce to the smooth, finite type case.

Firstly, if $X$ is any regular (not necessarily of finite type) $k$-scheme and $n\ge0$, then the sequence of pro \'etale sheaves
\begin{equation}
0\To\{W_r\Omega^n_{X,\sub{log}}\}_r\To\{W_r\Omega_X^n\}_r\xto{1-F}\{W_r\Omega_X^n\}_r\To0
\label{Illusie_fundamental}
\end{equation}
is known to be short exact by \cite[Cor.~2.9]{Shiho2007}, where $W_r\Omega^n_{X,\sub{log}}$ denotes the \'etale subsheaf of $W_r\Omega_X^n$ generated \'etale locally by logarithmic forms.

The continuous cohomology \cite{Jannsen1988} of a pro \'etale sheaf such as $\{W_r\Omega^n_{X,\sub{log}}\}_r$ will be denoted by $H^*_\sub{cont}(X_\sub{\'et},\{W_r\Omega^n_{X,\sub{log}}\}_r)$, or simply by $H^*_\sub{cont}(X_\sub{\'et},W\Omega^n_{X,\sub{log}})$ when there is no chance of confusion; similar notation is applied for continuous hypercohomology of pro complexes of \'etale sheaves, and in other topologies (if we do not specify a topology, it means Zariski). A more detailed discussion of such matters may be found in \cite[\S1.5.1]{Geisser2005}.

If $X$ is any $k$-scheme, then its crystalline cohomology groups $H^n_\sub{crys}(X/W)$ and $H^n_\sub{crys}(X):=H^n_\sub{crys}(X/W)\otimes_WK$ are defined using the crystalline site as in \cite{BerthelotOgus1978}; this does not require $X$ to satisfy any finite-type hypotheses. Illusie's comparison theorem \cite[\S2.I]{Illusie1979} states that there is a natural isomorphism \begin{equation}H^n_\sub{crys}(X/W)\isoto \bb H^n_\sub{cont}(X,W\Omega_X^\blob)\label{Illusie_comparision}\end{equation} for any smooth variety $X$ over $k$, but this remains true for any regular $k$-scheme $X$. (Proof: It is enough to show that $H^n_\sub{crys}(X/W_r)\isoto \bb H^n(X,W_r\Omega_X^\blob)$ for all $n\ge0$, $r\ge1$; by the Mayer--Vietoris property of both sides we may assume $X=\Spec C$ is affine. By N\'eron--Popescu desingularisation \cite{Popescu1985, Popescu1986} we may write $C=\indlim_\al C_\al$ as a filtered colimit of smooth $k$-algebras; since the de Rham--Witt complex commutes with filtered colimits, it is now enough to prove that $\indlim_\al H^n_\sub{crys}(C_\al/W_r)\isoto H^n_\sub{crys}(C/W_r)$ for all $n\ge0$, $r\ge1$. This is easily seen to be true if we pick compatible representations $C_\al=P_\al/J_\al$ of $C_\al$ as quotients of polynomial algebras over $W_r(k)$ and compute the crystalline cohomology in the usual way via the integrable connection arising on the divided power envelope of $(W_r(k),pW_r(k),\gamma)\to (P_\al,J_\al)$.) In particular, all proofs in Section \ref{subsection_main_local_proofs} will use $\bb H^n_\sub{cont}(X,W\Omega_X^\blob)$, but we will identify it with $H^n_\sub{crys}(X/W)$ when stating our results.

\subsection{Proofs of Theorems \ref{theorem_3} and \ref{theorem_4} via topological cyclic homology}\label{subsection_main_local_proofs}
We begin with a preliminary version of Theorem \ref{theorem_4}, phrased in terms of M.~Gros'  \cite{Gros1985} logarithmic crystalline Chern character $ch^\sub{log}:K_0(Y)_\bb Q\to\bigoplus_{i\ge0} H^i_\sub{cont}(Y, W\Omega_{Y,\sub{log}}^i)_\bb Q$. This is the most fundamental result of the article, depending essentially on a recent continuity theorem in topological cyclic homology.

\begin{proposition}\label{proposition_preliminary_main}
Let $X$ be a smooth, proper scheme over $k[[t_1,\dots,t_m]]$, and let $z\in K_0(Y)[\tfrac1p]$. Then the following are equivalent:
\begin{enumerate}
\item $z$ lifts to $(\projlim_rK_0(Y_r))[\tfrac1p]$.
\item $ch^\sub{log}(z)\in \bigoplus_{i\ge0}H^i_\sub{cont}(Y_\sub{\'et}, W\Omega_{Y,\sub{log}}^i)_\bb Q$ lifts to $\bigoplus_{i\ge0}H^i_\sub{cont}(X_\sub{\'et}, W\Omega_{X,\sub{log}}^i)_\bb Q$.
\end{enumerate}
\end{proposition}
\begin{proof}
The proof is an application of a continuity theorem in topological cyclic homology due to the author and B.~Dundas; for a summary of topological cyclic homology and its notation, we refer the reader to, e.g., \cite{Geisser2005} or \cite{Morrow_Dundas}. Indeed, according to \cite[Thm.~5.8]{Morrow_Dundas}, the canonical map $\TC(X;p)\to\op{holim}_r\TC(Y_r;p)$ is a weak equivalence, where $\TC(-;p)$ denotes the $p$-typical topological cyclic homology spectrum of a scheme; since the homotopy fibre of the trace map $tr:K(-)\to\TC(-;p)$ is nilinvariant by McCarthy's theorem \cite{McCarthy1997} (or rather the scheme-theoretic version of McCarthy's theorem coming from Zariski descent \cite{GeisserHesselholt1999}), it follows that there is a resulting homotopy cartesian square of spectra
\[\xymatrix{
\holim_rK(Y_r) \ar[d]\ar[r]& K(Y)\ar[d]^{tr} \\
\TC(X;p) \ar[r] & \TC(Y;p)
}\]
In other words, there is a commutative diagram of homotopy groups with exact rows:
\[\xymatrix{
\cdots\ar[r]&\pi_n\holim_rK(Y_r) \ar[d]\ar[r]& K_n(Y)\ar[d]^{tr}\ar[r]&\pi_{n-1}\holim_rK(Y_r,Y)\ar[d]^\cong \ar[r]&\cdots\\
\cdots\ar[r]&\TC_n(X;p) \ar[r] & \TC_n(Y;p)\ar[r] & \TC_{n-1}(X,Y;p) \ar[r]&\cdots\\
}\]
Inverting $p$, and noting that the natural map $\pi_n\holim_rK(Y_r)\to K_n(Y)$ factors through the surjection $\pi_n\holim_rK(Y_r)\to\projlim_rK_n(Y_r)$, one deduces from the diagram that the following are equivalent for any $z\in K_n(Y)[\tfrac1p]$:
\begin{enumerate}
\item $z$ lifts to $(\projlim_rK_n(Y_r))[\tfrac1p]$.
\item $tr(z)$ lifts to $\TC_n(X;p)[\tfrac1p]=\TC_n(X;p)_\bb Q$.
\end{enumerate}
Moreover, since $Y$ and $X$ are regular $\bb F_p$-schemes, there are natural decompositions \[\TC_n(Y;p)_\bb Q=\bigoplus_iH^{i-n}_\sub{cont}(Y_\sub{\'et}, W\Omega_{Y,\sub{log}}^i)_\bb Q,\qquad\TC_n(X;p)_\bb Q=\bigoplus_iH^{i-n}_\sub{cont}(X_\sub{\'et}, W\Omega_{X,\sub{log}}^i)_\bb Q,\] by T.~Geisser and L.~Hesselholt \cite[Thm.~4.1.1]{GeisserHesselholt1999}.

Taking $n=0$, the proof of the theorem will be completed as soon as we show that the rationalised trace map of topological cyclic homology \[tr:K_0(Y)_\bb Q\To \TC_0(Y;p)_\bb Q=\bigoplus_iH^i_\sub{cont}(Y_\sub{\'et}, W\Omega_{Y,\sub{log}}^i)_\bb Q\] is equal to Gros' logarithmic crystalline Chern character. This reduces via the usual splitting principle to the case of a line bundle $L\in H^1(Y,\roi_{Y}^\times)$, where it follows from the fact that both $tr(L)$ and the log crystalline Chern class $c_1^\sub{log}(L)$ are induced by the dlog map $\roi_{Y}^\times\to W\Omega^1_{Y}$; see \cite[Lem.~4.2.3]{GeisserHesselholt1999} and \cite[\S I.2]{Gros1985} respectively.
\end{proof}

To transform Proposition \ref{proposition_preliminary_main} into Theorem \ref{theorem_4}, we must study the relationship between the cohomology of the logarithmic de Rham--Witt sheaves and the eigenspaces of Frobenius acting on crystalline cohomology. For smooth, proper varieties over $k$, this follows from the general theory of slopes, but we require such comparisons also for the smooth, proper scheme $X$ over $k[[t_1,\dots,t_m]]$. To be more precise, for any regular $k$-scheme $X$ we denote by \[\ep:W_r\Omega_{X,\sub{log}}^i[-i]\To W_r\Omega_X^\blob\] the canonical map of complexes obtained from the inclusion $W_r\Omega_{X,\sub{log}}^i\subseteq W_r\Omega_X^i$, and we will consider the induced map on continuous cohomology. Note that $W_r\Omega_X^i$ has the same cohomology in the Zariski and \'etale topologies, since it is a quasi-coherent sheaf on the scheme $W_r(X)$ (e.g., by \cite[Prop.~2.18]{Shiho2007}); so we obtain an induced map \[\ep:H_\sub{cont}^i(X_\sub{\'et},W\Omega_{X,\sub{log}}^i)\To \bb H^{2i}_\sub{cont}(X_\sub{\'et},W\Omega_X^\blob)=\bb H^{2i}_\sub{cont}(X,W\Omega_X^\blob)=H^{2i}_\sub{crys}(X/W),\] (the final equality follows from line (\ref{Illusie_comparision})), which we study rationally in the following proposition:

\begin{proposition}\label{proposition_key}
Let $X$ be a regular $k$-scheme and $i\ge0$; consider the above map rationally: \[\ep_\bb Q:H_\sub{cont}^i(X_\sub{\'et},W\Omega_{X,\sub{log}}^i)_\bb Q\To H^{2i}_\sub{crys}(X).\] Then:
\begin{enumerate}
\item The image of $\ep_\bb Q$ is $H^{2i}_\sub{crys}(X)^{\phi=p^i}$.
\item If $X$ is a smooth, proper variety over a finite or algebraically closed field $k$, then $\ep_\bb Q$ is injective.
\end{enumerate}
\end{proposition}
\begin{proof}
Let $W_r\Omega_X^{\ge i}$ and $W_r\Omega_X^{<i}$ denote the naive upwards and downwards truncations of $W_r\Omega_X^\blob$ at degree $i$. Recalling the well-known de Rham--Witt identities  $dF=pFd$ and $Vd=pdV$, we may define morphisms of complexes \[\cal F:W_r\Omega_X^{\ge i}\To W_{r-1}\Omega_X^{\ge i},\qquad\cal V:W_r\Omega_X^{<i}\To W_{r+1}\Omega_X^{<i}\] degree-wise as \[p^\lambda F:W_r\Omega_X^{i+\lambda}\To W_{r-1}\Omega_X^{i+\lambda},\qquad p^{i-1-\lambda}V:W_r\Omega_X^{i-1-\lambda}\To W_{r+1}\Omega_X^{i-1-\lambda}\] for all $\lambda\ge0$. We claim that the resulting morphisms of pro complexes of sheaves \[1-\cal F:\{W_r\Omega_X^{>i}\}_r\To\{W_r\Omega_X^{>i}\}_r,\qquad 1-\cal V:\{W_r\Omega_X^{<i}\}_r\To\{W_r\Omega_X^{<i}\}_r\] are isomorphisms.

To prove this claim, it is sufficient to show that $1-p^\lambda F:\{W_r\Omega_X^n\}_r\to\{W_r\Omega_X^n\}$ and $1-p^{\lambda-1}V:\{W_r\Omega_X^n\}_r\to\{W_r\Omega_X^n\}$ are isomorphisms of pro sheaves for all $\lambda\ge1$. This follows from the fact that $p^\lambda F$ and $p^{\lambda-1}V$ are contracting operators; more precisely, inverses are provided by the maps \[\sum_{i=0}^{r-1}R^i(p^\lambda F)^{r-1-i}:W_{2r-1}\Omega_X^n\to W_r\Omega_X^n,\qquad \sum_{i=0}^{r-1}(p^{\lambda-1}RV)^{r-1-i}:W_r\Omega_X^n\to W_r\Omega_X^n.\]

Next, as recalled at line (\ref{Illusie_fundamental}), the sequence of pro \'etale sheaves \begin{equation*}
0\To\{W_r\Omega^i_{X,\sub{log}}\}_r\To\{W_r\Omega_X^i\}_r\xto{1-F}\{W_r\Omega_X^i\}_r\To0
\end{equation*}
is short exact; combining this with our observation that $1-\cal F$ is an isomorphism of $\{W_r\Omega_X^{>i}\}$, we arrive at a short exact sequence of pro complexes of sheaves
\begin{equation}
0\To\{W_r\Omega^i_{X,\sub{log}}\}_r\To\{W_r\Omega_X^{\ge i}\}_r\xto{1-\cal F}\{W_r\Omega_X^{\ge i}\}_r\To0,
\label{Illusie_fundamental2}
\end{equation}
which is well-known for smooth varieties over a perfect field \cite[\S I.3.F]{Illusie1979}. Recalling that $p^nF=\phi$ on $W_r\Omega_X^n$ we see that $p^i\cal F=\phi$, and so we obtain from (\ref{Illusie_fundamental2}) an exact sequence of rationalised continuous cohomology groups
\begin{equation}
H^i_\sub{cont}(X_\sub{\'et},W\Omega_{X,\sub{log}}^i)_\bb Q\To \bb H^{2i}_\sub{cont}(X,W\Omega_X^{\ge i})_\bb Q\xto{p^i-\phi} \bb H^{2i}_\sub{cont}(X,W\Omega_X^{\ge i})_\bb Q.
\label{log}
\end{equation}
 This proves that the canonical map 
\begin{equation}
 H^i_\sub{cont}(X_\sub{\'et},W\Omega_{X,\sub{log}}^i)_\bb Q\To \bb H^{2i}_\sub{cont}(X,W\Omega_X^{\ge i})_\bb Q^{\phi=p^i}
\label{surj1}
\end{equation}
is surjective.

Next, the short exact sequence of pro complexes of sheaves \[0\To\{W_r\Omega_X^{\ge i}\}_r\To \{W_r\Omega_X^{\blob}\}_r\To \{W_r\Omega_X^{<i}\}_r\To0\] gives rise to an exact sequence of rationalised continuous cohomology 
\begin{equation}
\bb H^{2i-1}_\sub{cont}(X,W\Omega_X^{<i})_\bb Q\to \bb H^{2i}_\sub{cont}(X,W\Omega_X^{\ge i})_\bb Q\to H^{2i}_\sub{crys}(X)\to \bb H^{2i}_\sub{cont}(X,W\Omega_X^{<i})_\bb Q.
\label{les>i}
\end{equation}
Recalling that $VF=FV=p$ on $W_r\Omega_X^n$, we see that $\phi\cal V=\cal V\phi=p^i$ on $W_r\Omega_X^{<i}$; since we have shown that $1-\cal V$ is an automorphism of $\bb H^n_\sub{cont}(X,W\Omega_X^{<i})$ for all $n\ge0$, it follows that $p^i-\phi$ is an is automorphism of $\bb H^n_\sub{cont}(X,W\Omega_X^{<i})_\bb Q$ for all $n\ge0$, so that in particular $p^i-\phi$ is an automorphism of the outer terms of the exact sequence (\ref{les>i}). By elementary linear algebra, the map of eigenspaces \begin{equation}\bb H^{2i}_\sub{cont}(X,W\Omega_X^{\ge i})_\bb Q^{\phi=p^i}\To H^{2i}_\sub{crys}(X)^{\phi=p^i}\label{surj}\end{equation} is therefore surjective.

Composing our two surjections proves (i). Now assume that $X$ is a smooth, proper variety over the perfect field $k$. Then the map of line (\ref{surj}) is injective, hence an isomorphism, by degeneration modulo torsion of the slope spectral sequence \cite[\S II.3.A]{Illusie1979}. So, to prove (ii), we must show that map (\ref{surj1}) is injective whenever $k$ is finite or algebraically closed. We begin with some general observations. Firstly, in (\ref{les>i}) we may replace $2i$ by $2i-1$ and then apply the same argument as immediately above to deduce that the map
\begin{equation}
\bb H^{2i-1}_\sub{cont}(X,W\Omega_X^{\ge i})_\bb Q^{\phi=p^i}\To H_\sub{crys}^{2i-1}(X)^{\phi=p^i}
\label{2i-1_espace}
\end{equation}
 is an isomorphism. Secondly, continuing (\ref{log}) to the left as a long exact sequence, we see that the kernel of (\ref{surj1}) is isomorphic to the cokernel of 
\begin{equation} 
\bb H^{2i-1}_\sub{cont}(X,W\Omega_X^{\ge i})_\bb Q\xto{p^i-\phi} \bb H^{2i-1}_\sub{cont}(X,W\Omega_X^{\ge i})_\bb Q,
\label{2i-1_espace2}
\end{equation}
so we must show that this map is surjective.

Now suppose that $k=\bb F_q$ is a finite field. According to the crystalline consequences of the Weil conjectures over finite fields \cite{KatzMessing1974}, the eigenvalues of the relative Frobenius $\phi_q:x\mapsto x^q$ acting on the crystalline cohomology $H^n_\sub{crys}(X)$ are algebraic over $\bb Q$ and all have complex absolute value $q^{n/2}$. It follows that $p^i$ cannot be an eigenvalue for the action of $\phi$ on $H^{2i-1}_\sub{crys}(X)$; so the right, hence the left, side of (\ref{2i-1_espace}) vanishes, and so map (\ref{2i-1_espace2}) is injective. But (\ref{2i-1_espace2}) is a $\bb Q_p$-linear endomorphism of a finite-dimensional $\bb Q_p$-vector space, so it must also be surjective, as desired.

Secondly suppose that $k=k^\sub{alg}$. Then $V:=\bb H^{2i-1}_\sub{cont}(X,W\Omega_X^{\ge i})_\bb Q$, equipped with operator $\cal F$, is an $F$-isocrystal over $k$. Since $k$ is algebraically closed, it is well-known that $1-\cal F:V\to V$ is therefore surjective. Indeed, by the Dieudonn\'e--Manin slope decomposition, we may suppose that $V=V_{r/s}$ is purely of slope $r/s$, where $r/s\in\bb Q_{\ge0}$ is a fraction written in lowest terms, i.e., $V=K^r$ and $\cal F^r=p^s\phi^r$: if $r\neq0$ it easily follows that $1-\cal F:V\to V$ is an automorphism; and if $r=0$ then $V=K$ and $\cal F=\phi$, whence $1-\cal F:V\to V$ is surjective as $k$ is closed under Artin--Schreier extensions.
\end{proof}

We may now prove Theorem \ref{theorem_4}; recall that $ch=ch^\sub{crys}$ denotes the crystalline Chern character:

\begin{theorem}\label{main_infinitesimal_theorem}
Let $X$ be a smooth, proper scheme over $k[[t_1,\dots,t_m]]$, where $k$ is a finite or algebraically closed field of characteristic $p$, and let $z\in K_0(Y)[\tfrac1p]$. Then the following are equivalent:
\begin{enumerate}
\item $z$ lifts to $(\projlim_r K_0(Y_r))[\tfrac1p]$
\item $ch(z)\in \bigoplus_{i\ge0}H_\sub{crys}^{2i}(Y)$ lifts to $\bigoplus_{i\ge0}H^{2i}_\sub{crys}(X)^{\phi=p^i}$.
\end{enumerate}
\end{theorem}
\begin{proof}
The proof is a straightforward diagram chase combining Propositions \ref{proposition_preliminary_main} and \ref{proposition_key} using the following diagram, in which we have deliberately omitted an unnecessary arrow:
\[\xymatrix{
(\projlim_rK_0(Y_r))[\tfrac1p]\ar[r] & K_0(Y)[\tfrac1p] \ar[d]^{ch^\sub{log}}\\
\bigoplus_iH^i_\sub{cont}(X_\sub{\'et}, W\Omega_{X,\sub{log}}^i)_\bb Q \ar@{->>}[d]^{\ep_\bb Q}\ar[r] & \bigoplus_iH^i_\sub{cont}(Y_\sub{\'et}, W\Omega_{Y,\sub{log}}^i)_\bb Q \ar[d]_\cong^{\ep_\bb Q}\\
\bigoplus_i H^{2i}_\sub{crys}(X)^{\phi=p^i} \ar[r] & \bigoplus_i H_\sub{crys}^{2i}(Y)^{\phi=p^i}
}\]
Given $z\in K_0(Y)[\tfrac1p]$, Proposition \ref{proposition_preliminary_main} states that $z$ lifts to $(\projlim_rK_0(Y_r))[\tfrac1p]$ if and only if $ch^\sub{log}(z)$ lifts to $\bigoplus_iH^i_\sub{cont}(X_\sub{\'et}, W\Omega_{X,\sub{log}}^i)_\bb Q$. However, Proposition \ref{proposition_key} states that the bottom left (resp.~right) vertical arrow is a surjection (resp.~an isomorphism); so the latter lifting condition is equivalent to $\ep_\bb Q(ch^\sub{log}(z))=ch(z)$ lifting to $\bigoplus_i H_\sub{crys}^{2i}(X)^{\phi=p^i}$, as required.
\end{proof}

In the case of a line bundle, Grothendieck's algebrization theorem allows us to prove a stronger result, thereby establishing Theorem \ref{theorem_3}:

\begin{corollary}\label{corollary_line_bundles}
Let $X$ be a smooth, proper scheme over $k[[t_1,\dots,t_m]]$, where $k$ is a finite or algebraically closed field of characteristic $p$, and let $L\in\Pic(Y)[\tfrac1p]$. Then the following are equivalent:
\begin{enumerate}
\item There exists $\tilde L\in\Pic(X)[\tfrac1p]$ such that $\tilde L|_Y=L$.
\item $c_1(L)\in H^2_\sub{crys}(Y)$ lifts to $H^2_\sub{crys}(X)^{\phi=p}$.
\end{enumerate}
\end{corollary}
\begin{proof}
The result follows from Theorem \ref{main_infinitesimal_theorem} and two observations: firstly, Grothendieck's algebrization theorem \cite[Thm.~5.1.4]{EGA_III_I} that $\Pic X=\projlim_r\Pic Y_r$; secondly, that $c_1(L)$ lifts to $H^2_\sub{crys}(X)^{\phi=p}$ if and only if $ch(L)=\exp(c_1(L))$ lifts to $\bigoplus_{i\ge0}H^{2i}_\sub{crys}(X)^{\phi=p^i}$.
\end{proof}

The final main result of this section is a modification of Corollary \ref{corollary_line_bundles} which extends Theorem \ref{theorem_VTC_for_divisors} to the base scheme $S=\Spec k[[t_1,\dots,t_m]]$; we do not provide all details of the proof:

\begin{theorem}\label{theorem_line_bundles}
Let $f:X\to S=\Spec k[[t_1,\dots,t_m]]$ be a smooth, projective morphism, where $k$ is a perfect field of characteristic $p>0$, let $L\in\Pic(Y)_\bb Q$, and let $c:=c_1(L)\in H^2_\sub{crys}(Y)$. Then the following are equivalent:
\begin{description}
\item[(deform)] There exists $\tilde L\in\Pic(X)_\bb Q$ such that $\tilde L|_Y=L$.
\item[(crys)] $c$ lifts to $H^2_\sub{crys}(X)$.
\item[(crys$\boldsymbol{-\phi}$)] $c$ lifts to $H^2_\sub{crys}(X)^{\phi=p}$.
\item[(flat)] $c$ is flat, i.e., it lifts to $H^0_\sub{crys}(S,R^2f_*\roi_{X/W})_\bb Q$.
\end{description}
\end{theorem}
\begin{proof}
We first claim that the analogue of Theorem \ref{theorem_Partie_Fixe} is true in this setting; that is, that the conditions (crys), (crys-$\phi$), and (flat) are equivalent for our local family $f:X\to \Spec k[[t_1,\dots,t_m]]$. The technical obstacle is that the theories of the convergent site and of isocrystals for non-finite-type schemes such as $X$ do not appear in the literature, though there is no doubt that the majority of these theories extend verbatim. To be precise, the key result we need is the following: if $d=\dim Y$ and $u\in H^2_\sub{crys}(X/W)$ denotes the Chern class of an ample line bundle, then the induced morphism of $\roi_{S/W}$-modules $u^i:R^{d-i}f_*\roi_{X/W}\to R^{d+i}f_*\roi_{X/W}$ is an isomorphism up to a bounded amount of $p$-torsion.

Using the arguments of \cite[Thm.~3.1]{Ogus1984} and \cite[Lem.~3.17]{Ogus1984} (which work in much greater generality than stated, since the base change theorem of crystalline cohomology \cite[Corol.~7.12]{BerthelotOgus1978} does not require the schemes to be of finite type over $k$), this isomorphism mod $p$-torsion may be checked on the special fibre $Y$. As in the proof of Theorem \ref{theorem_crystalline_partie_fixe}(i), this then follows from the Hard Lefschetz theorem for crystalline cohomology.

Deligne's axiomatic approach to the degeneration of Leray spectral sequences now shows that the rationalised Leray spectral sequence $E^{ab}_2=H^a_\sub{crys}(S,R^bf_*\roi_{X/W})_\bb Q\Rightarrow H^{a+b}_\sub{crys}(X)$ degenerates at the $E_2$-page, just as in the proof of Theorem \ref{theorem_crystalline_partie_fixe}(i). Verbatim repeating the rest of the proof of Theorem \ref{theorem_crystalline_partie_fixe}, and of Theorem \ref{theorem_Partie_Fixe}, shows that (crys), (crys-$\phi$), and (flat) are equivalent.

To complete the proof of the theorem, it remains to show that (crys-$\phi$) implies (deform); so assume that $c_1(L)$ lifts to $H^2_\sub{crys}(X)^{\phi=p}$. Write $A=k[[t_1,\dots,t_m]]$, whose strict Henselisation is $A^\sub{sh}=A\otimes_kk^\sub{alg}$, whose completion is $\hat{A^\sub{sh}}=k^\sub{alg}[[t_1,\dots,t_m]]$.

By applying Corollary \ref{corollary_line_bundles} and the same N\'eron--Popescu argument as in the proof of Theorem \ref{theorem_VTC_for_divisors} (whose indexing convention on line bundles we will follow), we find a line bundle $L_3\in\Pic(X\times_AA^\sub{sh})[\tfrac1p]$ whose restriction to $Y\times_kk^\sub{alg}$ coincides with the pullback of $L$ to $Y\times_kk^\sub{alg}$. Evidently there therefore exists a finite extension $k'$ of $k$ and a line bundle $L_5\in\Pic(X')[\tfrac1p]$, where $X':=X\times_Ak'[[t_1,\dots,t_m]]$, such that the restriction of $L_4$ to $Y\times_kk'$ coincides with the pullback of $L$ to $Y\times_kk'$.

Let $L_7\in\Pic(X)[\tfrac1p]$ be the pushforward of $L_4$ along the finite \'etale morphism $X'\to X$. Evidently $L_7|_Y=L^n$, where $n:= |k':k|$, and thus $\tilde L:=L_7^{1/n}\in\Pic(X)_\bb Q$ lifts $L$, as desired to prove (deform).
\end{proof}

\begin{remark}[Boundedness of $p$-torsion]\label{remark_boundedness}
The amount of $p$-torsion which must be neglected in Proposition \ref{proposition_preliminary_main} -- Corollary \ref{corollary_line_bundles} is bounded, i.e., annihilated by a large enough power of $p$ which depends only on $X$. For example, in Corollary \ref{corollary_line_bundles} there exists $\al\ge0$ with the following property: if $L\in\Pic(X)$ is such that $c_1(L)\in H^2_\sub{crys}(Y/W)$ lifts to $H^2_\sub{crys}(X/W)^{\phi=p}$, then $L^{p^\al}$ lifts to $\Pic(X)$.

To prove this boundedness claim, it is enough to observe that only a bounded amount of $p$-torsion must be neglected in both Proposition \ref{proposition_key} and the Geisser--Hesselholt decomposition of Proposition \ref{proposition_preliminary_main}. The former case is clear from the proof of Proposition \ref{proposition_key} and the finite generation of the $W$-modules $H^n_\sub{crys}(Y/W)$. The latter case is a consequence of the fact that the decomposition arises from a spectral sequence $E_2^{ab}=H^a_\sub{cont}(X_\sub{\'et},W\Omega_{X,\sub{log}}^{-b})\Rightarrow \TC_{-a-b}(X;p)$ (and similarly for $Y$); this spectral sequence is compatible with $\phi$, which acts as multiplication by $p^{-b}$ on $E_2^{ab}$, hence it degenerates modulo a bounded amount a $p$-torsion depending only on $\dim X$; see \cite[Thm.~4.1.1]{GeisserHesselholt1999}
\end{remark}

\begin{remark}[Lifting $L$ successively]\label{remark_step_bY_rtep}
Suppose that $X$ is a smooth, proper scheme over $k[[t_1,\dots,t_m]]$, and let $L\in\Pic(Y)$. Assuming that $L$ lifts to $L_r\in\Pic(Y_r)$ for some $r\ge1$, then there is a tautological obstruction in coherent cohomology to lifting $L_r$ to $\Pic(Y_{r+1})$, which lies in $H^2(Y,\roi_Y)$ if $m=1$. The naive approach to prove Theorem \ref{theorem_line_bundles} is to understand these tautological obstructions and thus successively lift $L$ to $\Pic(Y_2)$, $\Pic(Y_3)$, etc. (modulo a bounded amount of $p$-torsion). The proofs in this section have not used this approach, and in fact we strongly suspect that this naive approach does not work in general. We attempt to justify this suspicion in the rest of the remark.

If one unravels the details of the proof of Proposition \ref{proposition_preliminary_main} in the case of line bundles using the pro isomorphisms appearing in \cite{Morrow_Dundas}, one can prove the following modification of Corollary \ref{corollary_line_bundles}:

\begin{quote}
Given $r\ge1$ there exists $s\ge r$ such that for any $L\in\Pic(Y)$ the following implications hold (mod a bounded amount of $p$-torsion independent of~$r,s,L$):

$L$ lifts to $\Pic(Y_s)$$\implies$$c_1(L)$ lifts to $H^2_\sub{crys}(Y_s/W)^{\phi=p}$$\implies$$L$ lifts to $\Pic(Y_r)$
\end{quote}
We stress that $s$ is typically strictly bigger than $r$. In particular, it appears to be impossible to give a condition on $c_1(L)$, defined only in terms of $Y_r$, which ensures that $L$ lifts to $\Pic(Y_r)$.

Next, let us assume that $L$ is known to lift to $\tilde L\in\Pic(X)$ (we continue to ignore a bounded amount of $p$-torsion). Then the tautological coherent obstruction to lifting $L$ to $\Pic(Y_2)$ must vanish, so we may choose a lift $L_2\in\Pic(Y_2)$. Then it is entirely possible that $L_2\neq\tilde L|_{Y_2}$ and that $L_2$ does not lift to $\Pic(Y_3)$ (more precisely, this phenomenon would first occur at $Y_r$, for some $r\ge1$ depending on the amount of $p$-torsion neglected). However, see Remark \ref{remark_deJong}.

For these reasons it appears to be essential to prove the main results of this section by considering all infinitesimal thickenings of $Y$ at once, not one at a time. 
\end{remark}

\begin{remark}[de Jong's result]\label{remark_deJong}
The main implication (crys)$\Rightarrow$(deform) of Theorem \ref{theorem_constant_families} was proved by de Jong \cite{deJong2011} for smooth, proper $X$ over $k[[t]]$ under the following conditions: assumption (\ref{equ_assumption}) in Section \ref{Artin} holds, and $H^1(Y,\roi_Y)=H^0(Y,\Omega_{Y/k}^1)=0$.

In spite of Remark \ref{remark_step_bY_rtep}, de Jong proved his result by lifting the line bundle $L$ successively to $\Pic(Y_2)$, $\Pic(Y_3)$, etc. Indeed, the vanishing assumption $H^1(Y,\roi_Y)=0$ implies that all of the restriction maps \[\Pic(X)\to\cdots\to\Pic(Y_3)\to\Pic(Y_2)\to\Pic(Y)\] are injective, and therefore the problem discussed in the penultimate paragraph of Remark \ref{remark_step_bY_rtep} cannot occur: the arbitrarily chosen lift $L_2$ of $L$ must equal $\tilde L|_{Y_2}$, and hence $L_2$ lifts to $\Pic(Y_2)$, etc.
\end{remark}

\subsection{Cohomologically flat families over $k[[t]]$ and Artin's theorem}\label{Artin}
The aim of this section is to further analyse Corollary \ref{corollary_line_bundles} in a special case, to show that there is no obstruction to lifting line bundles in ``cohomologically constant'' families. This application of Corollary \ref{corollary_line_bundles} was inspired by de Jong's aforementioned work \cite{deJong2011}.

Let $k$ be an algebraically closed field of characteristic $p$, and let $f:X\to S=\Spec k[[t]]$ be a smooth, proper morphism.  Throughout this section we impose the following additional assumption on $X$, for which examples are given in Example \ref{example_of_free_crys_coh}:
\begin{equation}\mbox{The coherent $\roi_{S/W}$-modules $R^nf_*\roi_{X/W}$ are locally free, for all $n\ge0$.}\label{equ_assumption}\end{equation}
The assumption (and base change for crystalline cohomology; see especially \cite[Rmk.~7.10]{BerthelotOgus1978}) implies that each $\roi_{S/W}$-module $R^nf_*\roi_{X/W}$ is in fact an $F$-crystal on the crystalline site $(S/W)_\sub{crys}$; henceforth we simply say ``$F$-crystal over $k[[t]]$''.

We briefly review some of the theory of $F$-crystals. Let $\sigma:W[[t]]\to W[[t]]$ be the obvious lifting of the Frobenius on $k[[t]]$ which satisfies $\sigma(t)=t^p$. In the usual way \cite[\S2.4]{Katz1979}, we identify the category of $F$-crystals over $k[[t]]$ with the category of free, finite rank $W[[t]]$-modules $M$ equipped with a connection $\nabla:M\to M\,dt$ and a compatible $k[[t]]$-linear isogeny $F:\sigma^*M\to M$. The crystalline cohomology of the $F$-crystal is then equal to $H^*(M\xto{\nabla}M\,dt)$. An $F$-crystal $M=(M,F,\nabla)$ is {\em constant} if and only if it has the form $(\res M\otimes_WW[[t]],\res F\otimes\sigma,\tfrac{d}{dt})$ for some $F$-crystal $(\res M,\res F)$ over $k$.

We require the following simple lemma:

\begin{lemma}\label{lemma_Newton1}
Let $(M,F,\nabla)$ be an $F$-crystal over $k[[t]]$. Then:
\begin{enumerate}
\item The operator $p-F$ is an isogeny of $M\,dt$.
\item If $(M,F,\nabla)$ is constant then the composition $\ker\nabla\to M\to M/tM$ is an isomorphism of $W$-modules.
\end{enumerate}
\end{lemma}
\begin{proof}
(i): Since $F(m\,dt)=pt^{p-1}F(m)\,dt$, the operator $\tfrac1pF$ is well-defined on $M\,dt$ and is contracting. Therefore $1-\tfrac1pF$ is an automorphism of $M\,dt$, and so $p-F$ is injective with image $pM\,dt$.

(ii): Write $(M,F,\nabla)=(\res M\otimes_WW[[t]],\res F\otimes\sigma,\tfrac{d}{dt})$ for some $F$-crystal $(\res M,\res F)$ over $k$. Since $W[[t]]=\ker\tfrac{d}{dt}\oplus tW[[t]]$, it is easy to see that $M=\ker\nabla\oplus tM$, as desired.
\end{proof}

The following is the main theorem of the section, in which $k$ continues to be an algebraically closed field of characteristic $p$:

\begin{theorem}\label{theorem_constant_families}
Let $f:X\to S=\Spec k[[t]]$ be a smooth, proper morphism satisfying (\ref{equ_assumption}), and assume that the $F$-crystal $R^2f_*\roi_{X/W}$ on $k[[t]]$ is isogenous to a constant $F$-crystal. Then the cokernel of the restriction map $\Pic(X)\to\Pic(Y)$ is killed by a power of $p$.
\end{theorem}
\begin{proof}
We will prove this directly from Corollary \ref{corollary_line_bundles}, deliberately avoiding use of the stronger Theorem \ref{theorem_line_bundles}, at the expense of slightly lengthening the proof. The first half of the proof is similar to that of \cite[Thm.~1]{deJong2011}.

By assumption $R^2f_*\roi_{X/W}$ is isogenous to a constant $F$-crystal $(M,F,\nabla)$, and Lemma \ref{lemma_Newton1}(i) implies that the canonical map $\ker\nabla\to M/tM$ is an isomorphism. But, up to isogeny, the left side of this isomorphism is $H^0_\sub{crys}(S,R^2f_*\roi_{X/W})$ and the right side is $H^2_\sub{crys}(Y/W)$. In conclusion, the canonical map \[H^0_\sub{crys}(S,R^2f_*\roi_{X/W})_\bb Q\To H^2_\sub{crys}(Y)\] is an isogeny, and hence induces an isogeny on eigenspaces of the Frobenius.

Next note that $H^i_\sub{crys}(S,R^nf_*\roi_{X/W})=0$ for all $n\ge0$ and $i\ge 2$, since the crystalline cohomology of the crystal $R^nf_*\roi_{X/W}$ is computed using a two-term de Rham complex, as mentioned in the above review. Therefore the Leray spectral sequence for $f$ degenerates to short exact sequences \[0\To H^1_\sub{crys}(S,R^1f_*\roi_{X/W})\To H^2_\sub{crys}(X/W)\To H^0_\sub{crys}(S,R^2f_*\roi_{X/W})\To 0.\]

Letting $(M',F',\nabla')$ be the $F$-crystal $R^1f_*\roi_{X/W}$, Lemma \ref{lemma_Newton1} implies that the operator $p-F'$ is an isogeny of $M\,dt$, hence is surjective modulo a bounded amount of $p$-torsion on its quotient $\op{Coker}(M\xto{\nabla}M\,dt)$, which is isogenous to $H^1_\sub{crys}(S,R^1f_*\roi_{X/W})_\bb Q$. It follows from elementary linear algebra that the surjection $H^2_\sub{crys}(X/W)\to H^0_\sub{crys}(S,R^2f_*\roi_{X/W})$ remains surjective modulo a bounded amount of $p$-torsion after restricting to $p$-eigenspaces of the Frobenius.

Combining the established surjection and isogeny proves that the canonical map $H^2_\sub{crys}(X/W)^{\phi=p}\to H^2_\sub{crys}(Y/W)^{\phi=p}$ is surjective modulo a bounded amount of $p$-torsion. Corollary \ref{corollary_line_bundles}, or rather its improvement explained in Remark \ref{remark_boundedness}, completes the proof.
\end{proof}

\begin{remark}\label{remark_NS_p_torsion}
Let $f:X\to S=\Spec k[[t]]$ satisfy the assumptions of Theorem \ref{theorem_constant_families}; then a consequence of the theorem is that the cokernel of the N\'eron--Severi specialisation map \cite[Exp.~X, \S7]{SGA_VI} \[sp:\NS(X\times_Sk((t))^\sub{alg})\To\NS(Y)\] is a finite $p$-group. Indeed, $sp$ arises as a quotient of the colimit of the maps \begin{equation}\Pic(X\times_S F)\cong\Pic(X\times_SA)\To\Pic(Y),\label{eqn_pic}\end{equation} where $F$ varies over all finite extensions of $k((t))$ inside $k((t))^\sub{alg}$ and $A$ is the integral closure of $k[[t]]$ inside $F$. But the hypotheses of Theorem \ref{theorem_constant_families} are satisfied for each morphism $X\times_SA\to\Spec A$, by base change for crystalline cohomology, and hence the cokernel of (\ref{eqn_pic}) is killed by a power of $p$. Passing to the colimit over $F$ we deduce that the cokernel of $sp$ is a $p$-torsion group; it is moreover finite since $\NS(Y)$ is finitely generated.
\end{remark}

\begin{example}\label{example_of_free_crys_coh}
Let $f:X\to S=\Spec k[[t]]$ be a smooth, proper morphism. Then assumption (\ref{equ_assumption}) is known to be satisfied in each of the following cases:
\begin{enumerate}\itemsep0pt
\item $X$ is a family of K3 surfaces over $S$.
\item $X$ is an abelian scheme over $S$ \cite[Corol.~2.5.5]{BerthelotBreenMessing1982}.
\item $X$ is a complete intersection in $\bb P_S^d$. (Proof: Given any affine pd-thickening $S'$ of $S$, there exists a lifting of $X$ to a smooth complete intersection $X'$ in $\bb P^d_{S'}$; then apply Berthelot's comparison isomorphism $R^nf_*\roi_{X/W}(S')\cong H^n_\sub{dR}(X'/S')$, which is locally free since $X'$ is a smooth complete intersection over $S'$ \cite[Exp.~XI, Thm.~1.5]{SGA_VII_I}.)
\end{enumerate}
Moreover, under assumption (\ref{equ_assumption}), $R^2f_*\roi_{X/W}$ is isogenous to a constant $F$-crystal if the geometric generic fibre of $X$ is supersingular in degree $2$, i.e., the crystalline cohomology $H^2_\sub{crys}(X\times_{k[[t]]}k((t))^\sub{perf}/W(k((t))^\sub{perf}))$ is purely of slope $1$. Indeed, generic supersingularity forces the Newton polygon of $R^2f_*\roi_{X/W}$ to be purely of slope $1$ at both geometric points of $S$, and so $R^2f_*\roi_{X/W}$ is isogenous to a constant $F$-crystal by \cite[Thm.~2.7.1]{Katz1979}.

In particular, Theorem \ref{theorem_constant_families} and Remark \ref{remark_NS_p_torsion} apply to any smooth, proper family of K3 surfaces over $k[[t]]$ whose geometric generic fibre is supersingular; this reproves a celebrated theorem of M.~Artin \cite[Corol.~1.3]{Artin1974}.
\end{example}

\section{An application to the Tate conjecture}\label{section_application}
In this section we apply Theorem \ref{theorem_2} to the study of the Tate conjecture; more precisely, we reduce the Tate conjecture for divisors to the case of surfaces. This result is already known to some experts, and our proof is a modification of de Jong's \cite{deJong20??}. Throughout this section $k$ denotes a finite field. We begin with the following folklore result that all formulations of this conjecture are equivalent:

\begin{proposition}\label{proposition_Tate_conj}
Let $X$ be a smooth, projective, geometrically connected variety over a finite field $k=\bb F_{p^m}$, and let $\ell\neq p$ be a prime number. Then the following are equivalent:
\begin{enumerate}
\item The $\ell$-adic Chern class $c_1^\sub{\'et}\otimes\bb Q_\ell:\Pic(X)_{\bb Q_\ell}\to H^2_\sub{\'et}(X\times_kk^\sub{alg},\bb Q_\ell(1))^{\Gal(k^\sub{alg}/k)}$ is surjective.
\item The crystalline Chern class $c_1^\sub{crys}\otimes K:\Pic(X)_{K}\to H^2_\sub{crys}(X)^{\phi^m=p^m}$ is surjective.
\item The crystalline Chern class $c_1^\sub{crys}\otimes\bb Q_p:\Pic(X)_{\bb Q_p}\to H^2_\sub{crys}(X)^{\phi=p}$ is surjective.
\item The order of the pole of the zeta function $\zeta(X,s)$ at $s=1$ is equal to the dimension of $A^1_\sub{num}(X)_\bb Q$, the group of rational divisors modulo numerical equivalence.
\end{enumerate}
\end{proposition}
\begin{proof}
We start with some remarks about the Frobenius in crystalline and \'etale cohomology; let $q:=p^m$. Firstly, $\phi^m:x\mapsto x^q$ is an endomorphism of the $k$-variety $X$, which induces an endomorphism $\phi^m\otimes1$ of the $k^\sub{alg}$-variety $\res X:=X\times_kk^\sub{alg}$, which in turn induces a $\bb Q_\ell$-linear automorphism $F$ of $H^n_\sub{\'et}(\res X,\bb Q_\ell)$ for any $n\ge0$ and any prime $\ell\neq p$; it is well-known (see, e.g., \cite[Rmk.~13.5]{Milne1980}) that $F$ is the same as the automorphism induced by $1\otimes\sigma$, where $\sigma\in\Gal(k^\sub{alg}/k)$ is the geometric Frobenius relative to $k$. Fixing for convenience an isomorphism $\bb Z_\ell(1)\cong\bb Z_\ell$, it follows that the $\ell$-adic cycle class map $cl^{(\ell)}_i:\CH^i(X)_\bb Q\to H^{2i}_\sub{\'et}(\res X,\bb Q_\ell)$ has image inside the eigenspace $H^{2i}_\sub{\'et}(\res X,\bb Q_\ell)^{F=q^i}\cong H^{2i}_\sub{\'et}(\res X,\bb Q_\ell(i))^{\Gal(k^\sub{alg}/k)}$, and the surjectivity part of the $\ell$-adic Tate conjecture is the assertion that its image spans the entire eigenspace.

Meanwhile in crystalline cohomology, the $\bb Q_p$-linear automorphism $\phi$ of $H^*_\sub{crys}(X)$ induces the $K$-linear automorphism $F:=\phi^m$ (hopefully using the same $F$ to denote the Frobenius in both crystalline and $\ell$-adic \'etale cohomology will not be a source of confusion).

To summarise, let $\ell$ denote any prime number (possibly equal to $p$), let \[cl_i^{(\ell)}:\CH^i(X)_\bb Q\To H^{2i}_{(\ell)}(X):=\begin{cases} H^{2i}_\sub{\'et}(\res X,\bb Q_\ell) & \ell\neq p \\ H^{2i}_\sub{crys}(X) & \ell=p\end{cases}\] denote the associated cycle class map (whose image lies in $H^{2i}_{(\ell)}(X)^{F=q^i}$), and let \[\Lambda:=\begin{cases}\bb Q_\ell & \ell\neq p \\ K& \ell=p\end{cases}\] be the coefficient field of the Weil cohomology theory. The Chow group of rational cycles modulo homological equivalence is \[A^i_{\ell\sub{-hom}}(X)_\bb Q:=\CH^i(X)_\bb Q/\ker cl_i^{(\ell)},\] and we denote by $cl_i^\Lambda:A^i_{\ell\sub{-hom}}(X)_\bb Q\otimes_\bb Q\Lambda\to H^{2i}_{(\ell)}(X)^{F=q^i}$ the induced cycle class map. We may finally uniformly formulate statements (i) and (ii) as the assertion that $cl_1^\Lambda$ is surjective.

We may now properly begin the proof, closely follow arguments from \cite{Tate1994}, though Tate considered only the case of \'etale cohomology. Homological equivalence (for either $\ell$-adic \'etale cohomology, $\ell\neq p$, or crystalline cohomology) and numerical equivalence for divisors agree rationally \cite[Prop.~3.4.6.1]{Andre2004}; this implies that $A^1_{\ell\sub{-hom}}(X)_\bb Q$ is independent of the prime $\ell$ and that the intersection pairing \[A^1_{\ell\sub{-hom}}(X)_\bb Q\times A^{d-1}_{\ell\sub{-hom}}(X)_\bb Q\To A^d_{\ell\sub{-hom}}(X)_\bb Q\xto{\sub{deg}}\bb Q\] is non-degenerate on the left. Combining these observations, we obtain a commutative diagram

\[\xymatrix@C=5mm@R=10mm{
H^2_{(\ell)}(X)^{F=q} \ar@{^(->}[r]^j & H^2_{(\ell)}(X)^{F=q}_\sub{gen} \ar[r]^\gamma & H^2_{(\ell)}(X)_{F=q}\ar[r]^-\cong & \Hom_\Lambda(H^{2d-2}_{(\ell)}(X)^{F=q},\Lambda)\ar[d]^{\sub{dual of }cl_{d-1}^{\Lambda}}\\
A^1_{\ell\sub{-hom}}(X)_\bb Q\otimes_\bb Q\Lambda\ar[u]^{cl_1^\Lambda} \ar@{^(->}[r]&\mathrlap{\Hom_\bb Q(A^{d-1}_{\ell\sub{-hom}}(X)_\bb Q,\bb Q)\otimes_\bb Q\Lambda}\rule{15mm}{0pt} &\rule{26mm}{0pt}\ar@{^(->}[r]& \Hom_\Lambda(A^{d-1}_{\ell\sub{-hom}}(X)_\bb Q\otimes_\bb Q\Lambda,\Lambda)
}\]
where:
\begin{itemize}\itemsep0pt
\item The top right (resp.~bottom left) horizontal arrow arises from the pairing on cohomology groups (resp.~Chow groups); the top right horizontal arrow is an isomorphism by the compatibility of Poincar\'e duality with the action of $F$.
\item $H^2_{(\ell)}(X)^{F=q}_\sub{gen}\subseteq H^2_{(\ell)}(X)$ is the generalised eigenspace for $F$ of eigenvalue $q$.
\item $H^2_{(\ell)}(X)_{F=q}:=H^2_{(\ell)}(X)/\Im(q-F)$.
\item $\gamma$ is the composition $H^2_{(\ell)}(X)^{F=q}_\sub{gen}\into H^2_{(\ell)}(X)\onto H^2_{(\ell)}(X)_{F=q}$.
\end{itemize}
It follows from the diagram that $cl^\Lambda_1$ is injective. The proof can now be quickly completed.

According to the Lefschetz trace formula (for either $\ell$-adic \'etale cohomology where $\ell\neq p$, or crystalline cohomology \cite[Thm.~1]{KatzMessing1974}), the order of the pole of interest in (iv) is equal to $\dim_\Lambda H^2_{(\ell)}(X)^{F=q}_\sub{gen}$. The equality $A^1_{\ell\sub{-hom}}(X)_\bb Q=A^1_\sub{num}(X)_\bb Q$ and the injectivity of $cl^\Lambda_1$ and $j$ therefore immediately imply: \begin{center}Statement (iv) is true $\iff$ $cl^\Lambda_1$ and $j$ are surjective.\end{center}
Moreover, by elementary linear algebra and diagram chasing:
\begin{center}
\begin{tabular}{lcl}
$cl^\Lambda_1$ is surjective & $\iff$ & $cl^\Lambda_1$ is an isomorphism\\
&$\implies$& $\gamma j$ is injective, i.e.,~$\ker(q-F)\cap\Im(q-F)=0$\\
&$\iff$&$j$ is surjective, i.e.,~$\ker(q-F)=\ker(q-F)^N$ for $N\ge1$
\end{tabular}
\end{center}
In conclusion, statement (iv) is equivalent to the surjectivity of $cl^\Lambda_1$, as required to prove (i)$\Leftrightarrow$(iv) and (ii)$\Leftrightarrow$(iv).

Finally, the equivalence of (ii) and (iii) is an immediate consequence of Galois descent for vector spaces, which implies that the canonical map \[H^2_\sub{crys}(X)^{\phi=p}\otimes_{\bb Q_p}K\To H^2_\sub{crys}(X)^{F=q}\] is an isomorphism.
\end{proof}

We will refer to the equivalent statements of Proposition \ref{proposition_Tate_conj} as ``the Tate conjecture for divisors for $X$'' or, following standard terminology, simply as $T^1(X)$. The main theorem of this section is the reduction of the Tate conjecture for divisors to the case of surfaces. We begin with a discussion of the problem. If $X$ is a smooth, projective, geometrically connected $k$-variety of dimension $d>3$ and $Y\into X$ is a smooth ample divisor (which always exists after base changing by a finite extension $k'$ of $k$; note that $T^1(X\times_kk')\Rightarrow T^1(X)$), then the restriction maps $H^2_\sub{crys}(X)\to H^2_\sub{crys}(Y)$ and $\Pic(X)_\bb Q\to\Pic(Y)_\bb Q$ are isomorphisms (see Remark \ref{remark_hyerplane_theorem}), and so it follows immediately that $T^1(X)$ and $T^1(Y)$ are equivalent. This reduces the Tate conjecture for divisors to $3$-folds and surfaces.

Now suppose $\dim X=3$ and continue to let $Y\into X$ be a smooth ample divisor. If Corollary \ref{corollary_hypersurfaces} were known to be true with $\Pic(Y)_\bb Q$ and $\Pic(X)_\bb Q$ replaced by $\Pic(Y)_{\bb Q_p}$ and $\Pic(X)_{\bb Q_p}$, then $T^1(X)$ would easily follow from $T^1(Y)$, thereby proving the desired reduction to surfaces. Unfortunately, the analogue of Corollary \ref{corollary_hypersurfaces} with $\bb Q_p$-coefficients is currently unknown. Instead, we will choose the smooth ample divisor $Y$ so that the $p$-eigenspace of its crystalline cohomology as small as possible; we borrow this style of argument from the aforementioned earlier proof of Theorem \ref{theorem_reduction_of_Tate} by de Jong, who used deep monodromy results proved with N.~Katz \cite{deJongKatz2000} to choose $Y$ more stringently:

\begin{lemma}\label{lemma_Chebotaryev}
Let $X$ be a smooth, projective, geometrically connected, $3$-dimensional variety over a finite field $k=\bb F_{p^m}$. Then there exists a finite extension $k'=\bb F_{p^M}$ of $k$ and a smooth ample divisor $Y\into X\times_kk'$ such that the following inequality of dimensions of generalised eigenspaces holds: \[\dim_KH^2_\sub{crys}(X)^{\phi^m=p^m}_\sub{gen}\ge\dim_{K'}H^2_\sub{crys}(Y)^{\phi^M=p^M}_\sub{gen}\] (where $K':=\Frac W(k')$).
\end{lemma}
\begin{proof}
Let $\{X_t:t\in\bb P^1_k\}$ be a pencil of hypersurface sections of degree $\ge 2$, and recall the crystalline version \cite{KatzMessing1974} of Deligne's ``Least common multiple theorem'' \cite[Thm.~4.5.1]{Deligne1980}, which states that the polynomial $\det(1-\phi^mT|H^2_\sub{crys}(X))$ is equal to the least common multiple of all complex polynomials $\prod_i(1-\al_iT)\in 1+T\bb C[T]$ with the following property: for all $M\in m\bb Z$ and all $t\in\bb P^1_k(\bb F_{p^M})$ such that $X_t$ is smooth, the polynomial $\prod_i(1-~\al^{M/m}T)$ divides $\det(1-\phi^MT|H^2_\sub{crys}(X_t))$.

It follows that there exists $M\in m\bb Z$ and $t\in\bb P_k^1(\bb F_{p^M})$ such that $X_t$ is smooth and such that the multiplicity of $1-p^MT$ as a factor of $\det(1-\phi^MT|H^2_\sub{crys}(X_t))$ is at most the multiplicity of $1-p^mT$ as a factor of $\det(1-\phi^mT|H^2_\sub{crys}(X))$. But these multiplicities are precisely the dimensions of the generalised eigenspaces of interest.
\end{proof}

\begin{theorem}\label{theorem_reduction_of_Tate}
Let $k$ be a finite field, and assume that $T^1(X)$ is true for every smooth, projective, connected surface over $k$. Then $T^1(X)$ is true for every smooth, projective, connected variety over $k$.
\end{theorem}
\begin{proof}
As we discussed above, it is sufficient to show that $T^1$ for surfaces implies $T^1$ for $3$-folds. So let $X$ be a smooth, projective, connected $3$-fold over $k=\bb F_{p^m}$; replacing $k$ by $H^0(X,\roi_X)$, we may assume that $X$ is geometrically connected. Let $k'=\bb F_{p^M}$ be a finite extension of $k$ and  $Y\into X':=X\times_kk'$ a smooth ample divisor as in the statement of Lemma \ref{lemma_Chebotaryev}.

Consider the following diagram of inclusions of $K':=\Frac W(k')$ vector spaces:
\[\xymatrix@=4mm{
H^2_\sub{crys}(X)^{\phi^m=p^m}_\sub{gen}\otimes_KK' \ar@{}[r]|-{\mbox{$\subseteq$}}& H^2_\sub{crys}(X')^{\phi^M=p^M}_\sub{gen} \ar@{}[r]|-{\mbox{$\subseteq$}}& H^2_\sub{crys}(Y)^{\phi^M=p^M}_\sub{gen}\\
H^2_\sub{crys}(X)^{\phi^m=p^m}\otimes_KK' \ar@{}[u]|-{\rotatebox{90}{$\subseteq$}} \ar@{}[r]|-{\mbox{$\subseteq$}}& H^2_\sub{crys}(X')^{\phi^M=p^M} \ar@{}[u]|-{\rotatebox{90}{$\subseteq$}}\ar@{}[r]|-{\mbox{$\subseteq$}}& H^2_\sub{crys}(Y)^{\phi^M=p^M}\ar@{}[u]|-{\rotatebox{90}{$\subseteq$}}
}\]
The vertical inclusions are obvious; the left horizontal inclusions result from the identification $H^2_\sub{crys}(X')=H^2_\sub{crys}(X)\otimes_KK'$; the right horizontal inclusions are a consequence of Weak Lefschetz for crystalline cohomology \cite{KatzMessing1974}, stating that the restriction map $H^2_\sub{crys}(X')\to H^2_\sub{crys}(Y)$ is injective.

Since we are assuming the validity of $T^1(Y)$, it follows from the proof of Proposition \ref{proposition_Tate_conj} that $H^2_\sub{crys}(Y)^{\phi^M=p^M}_\sub{gen}=H^2_\sub{crys}(Y)^{\phi^M=p^M}$; hence $H^2_\sub{crys}(X')^{\phi^M=p^M}_\sub{gen}=H^2_\sub{crys}(X')^{\phi^M=p^M}$. Moreover, by choice of $Y$, the dimension of the top-right entry in the diagram is at most the dimension of the top-left entry. It follows that the inclusions of generalised eigenspaces at the top of the diagram are in fact equalities, and hence that $H^2_\sub{crys}(X')^{\phi^M=p^M}=H^2_\sub{crys}(Y)^{\phi^M=p^M}$. More precisely, we have proved that the restriction map $H^2_\sub{crys}(X')\to H^2_\sub{crys}(Y)$ induces an isomorphism $H^2_\sub{crys}(X')^{\phi^M=p^M}\isoto H^2_\sub{crys}(Y)^{\phi^M=p^M}$, from which it follows that \[H^2_\sub{crys}(X')^{\phi=p}\isoto H^2_\sub{crys}(Y)^{\phi=p}\] by the same Galois descent argument as at the end of the proof of Proposition \ref{proposition_Tate_conj}.

But now it is an easy consequence of Corollary \ref{corollary_hypersurfaces} that the restriction map \[\Im(\Pic(X)_\bb Q\xto{c_1} H^2_\sub{crys}(X)^{\phi=p})\To \Im(\Pic(Y)_\bb Q\xto{c_1} H^2_\sub{crys}(Y)^{\phi=p})\] is also an isomorphism. Since the $\bb Q_p$-linear span of the right image is $H^2_\sub{crys}(Y)^{\phi=p}$, by $T^1(Y)$, it follows that the $\bb Q_p$-linear span of the left image of $H^2_\sub{crys}(X')^{\phi=p}$; i.e., $T^1(X')$ is true. Finally recall that $T^1(X')\Rightarrow T^1(X)$, completing the proof.
\end{proof}


\begin{thebibliography}{10}

\bibitem{Andre2004}
{\sc {Andr\'e}, Y.}
\newblock {\em {Une introduction aux motifs. Motifs purs, motifs mixtes,
  p\'eriodes.}}
\newblock Paris: Soci\'et\'e Math\'ematique de France, 2004.

\bibitem{Artin1974}
{\sc Artin, M.}
\newblock Supersingular {$K3$} surfaces.
\newblock {\em Ann. Sci. \'Ecole Norm. Sup. (4) 7\/} (1974), 543--567 (1975).

\bibitem{Berthelot1986}
{\sc Berthelot, P.}
\newblock G\'eom\'etrie rigide et cohomologie des vari\'et\'es alg\'ebriques de
  caract\'eristique {$p$}.
\newblock {\em M\'em. Soc. Math. France (N.S.)}, 23 (1986), 3, 7--32.
\newblock Introductions aux cohomologies $p$-adiques (Luminy, 1984).

\bibitem{BerthelotBreenMessing1982}
{\sc Berthelot, P., Breen, L., and Messing, W.}
\newblock {\em Th\'eorie de {D}ieudonn\'e cristalline. {II}}, vol.~930 of {\em
  Lecture Notes in Mathematics}.
\newblock Springer-Verlag, Berlin, 1982.

\bibitem{SGA_VI}
{\sc Berthelot, P., Grothendieck, A., and Illusie, L.}
\newblock {\em Th\'eorie des intersections et th\'eor\`eme de
  {R}iemann-{R}och}.
\newblock Lecture Notes in Mathematics, Vol. 225. Springer-Verlag, Berlin-New
  York.
\newblock S{\'e}minaire de G{\'e}om{\'e}trie Alg{\'e}brique du Bois-Marie
  1966--1967 (SGA 6).

\bibitem{BerthelotIllusie1970}
{\sc Berthelot, P., and Illusie, L.}
\newblock Classes de {C}hern en cohomologie cristalline.
\newblock {\em C. R. Acad. Sci. Paris S\'er. A-B 270 (1970), A1695-A1697; ibid.
  270\/} (1970), A1750--A1752.

\bibitem{BerthelotOgus1978}
{\sc {Berthelot}, P., and {Ogus}, A.}
\newblock {Notes on crystalline cohomology.}
\newblock {Mathematical Notes. Princeton, New Jersey: Princeton University
  Press. Tokyo: University of Tokyo Press. VI, not consecutively paged. \$ 9.50
  (1978).}, 1978.

\bibitem{BerthelotOgus1983}
{\sc Berthelot, P., and Ogus, A.}
\newblock {$F$}-isocrystals and de {R}ham cohomology. {I}.
\newblock {\em Invent. Math. 72}, 2 (1983), 159--199.

\bibitem{BlochEsnaultKerz2012}
{\sc {Bloch}, S., {Esnault}, H., and {Kerz}, M.}
\newblock {$p$-adic deformation of algebraic cycle classes}.
\newblock {\em Invent. Math.\/} (2012).
\newblock Awaiting publication.

\bibitem{BlochEsnaultKerz2013}
{\sc {Bloch}, S., {Esnault}, H., and {Kerz}, M.}
\newblock {Deformation of algebraic cycle classes in characteristic zero}.
\newblock {\em Algebraic Geometry 1}, 3 (2014), 290--310.

\bibitem{LangerZinkDavis2011}
{\sc Davis, C., Langer, A., and Zink, T.}
\newblock Overconvergent de {R}ham-{W}itt cohomology.
\newblock {\em Ann. Sci. \'Ec. Norm. Sup\'er. (4) 44}, 2 (2011), 197--262.

\bibitem{deJong20??}
{\sc de~Jong, A.~J.}
\newblock Tate conjecture for divisors.
\newblock {\em Unpublished note\/}.

\bibitem{deJong1996}
{\sc de~Jong, A.~J.}
\newblock Smoothness, semi-stability and alterations.
\newblock {\em Inst. Hautes \'Etudes Sci. Publ. Math.}, 83 (1996), 51--93.

\bibitem{deJong2011}
{\sc de~Jong, A.~J.}
\newblock On a result of {A}rtin.
\newblock {\em Available at {\tt\url{www.math.columbia.edu/~dejong/}}\/}
  (2011).

\bibitem{deJongKatz2000}
{\sc de~Jong, A.~J., and Katz, N.~M.}
\newblock Monodromy and the {T}ate conjecture: {P}icard numbers and
  {M}ordell-{W}eil ranks in families.
\newblock {\em Israel J. Math. 120}, part A (2000), 47--79.

\bibitem{Deligne1968}
{\sc {Deligne}, P.}
\newblock {Th\'eor\`eme de Lefschetz et crit\`eres de d\'eg\'en\'erescence de
  suites spectrales.}
\newblock {\em {Publ. Math., Inst. Hautes \'Etud. Sci.} 35\/} (1968), 107--126.

\bibitem{Deligne1971}
{\sc Deligne, P.}
\newblock Th\'eorie de {H}odge. {II}.
\newblock {\em Inst. Hautes \'Etudes Sci. Publ. Math.}, 40 (1971), 5--57.

\bibitem{Deligne1980}
{\sc {Deligne}, P.}
\newblock {La conjecture de Weil. II.}
\newblock {\em {Publ. Math., Inst. Hautes \'Etud. Sci.} 52\/} (1980), 137--252.

\bibitem{Morrow_Dundas}
{\sc Dundas, B., and Morrow, M.}
\newblock Finite generation and continuity of topological {H}ochschild and
  cyclic homology.
\newblock {\em {\tt arXiv:1403.0534}\/} (2014).

\bibitem{Geisser2005}
{\sc Geisser, T.}
\newblock Motivic cohomology, {$K$}-theory and topological cyclic homology.
\newblock In {\em Handbook of {$K$}-theory. {V}ol. 1, 2}. Springer, Berlin,
  2005, pp.~193--234.

\bibitem{GeisserHesselholt1999}
{\sc {Geisser}, T., and {Hesselholt}, L.}
\newblock {Topological cyclic homology of schemes.}
\newblock In {\em {Algebraic $K$-theory. Proceedings of an AMS-IMS-SIAM summer
  research conference, Seattle, WA, USA, July 13--24, 1997}}. Providence, RI:
  American Mathematical Society, 1999, pp.~41--87.

\bibitem{GilletMessing1987}
{\sc Gillet, H., and Messing, W.}
\newblock Cycle classes and {R}iemann-{R}och for crystalline cohomology.
\newblock {\em Duke Math. J. 55}, 3 (1987), 501--538.

\bibitem{Gros1985}
{\sc Gros, M.}
\newblock Classes de {C}hern et classes de cycles en cohomologie de
  {H}odge-{W}itt logarithmique.
\newblock {\em M\'em. Soc. Math. France (N.S.)}, 21 (1985), 1--87.

\bibitem{EGA_III_I}
{\sc Grothendieck, A.}
\newblock \'{E}l\'ements de g\'eom\'etrie alg\'ebrique. {III}. \'{E}tude
  cohomologique des faisceaux coh\'erents. {I}.
\newblock {\em Inst. Hautes \'Etudes Sci. Publ. Math.}, 11 (1961), 5--167.

\bibitem{Grothendieck1966}
{\sc Grothendieck, A.}
\newblock On the de {R}ham cohomology of algebraic varieties.
\newblock {\em Publ. Math. Inst. Hautes \'Etudes Sci.}, 29 (1966), 95--103.

\bibitem{SGA_II}
{\sc Grothendieck, A.}
\newblock {\em Cohomologie locale des faisceaux coh\'erents et th\'eor\`emes de
  {L}efschetz locaux et globaux {$(SGA$} {$2)$}}.
\newblock North-Holland Publishing Co., Amsterdam; Masson \& Cie, \'Editeur,
  Paris, 1968.

\bibitem{SGA_VII_I}
{\sc Grothendieck, A.}
\newblock {\em Groupes de monodromie en g\'eom\'etrie alg\'ebrique. {I}}.
\newblock Lecture Notes in Mathematics, Vol. 288. Springer-Verlag, Berlin-New
  York, 1972.
\newblock S{\'e}minaire de G{\'e}om{\'e}trie Alg{\'e}brique du Bois-Marie
  1967--1969 (SGA 7 I).

\bibitem{Illusie1975}
{\sc Illusie, L.}
\newblock Report on crystalline cohomology.
\newblock In {\em Algebraic geometry ({P}roc. {S}ympos. {P}ure {M}ath., {V}ol.
  29, {H}umboldt {S}tate {U}niv., {A}rcata, {C}alif., 1974)}. Amer. Math. Soc.,
  Providence, R.I., 1975, pp.~459--478.

\bibitem{Illusie1979}
{\sc Illusie, L.}
\newblock Complexe de de\thinspace {R}ham-{W}itt et cohomologie cristalline.
\newblock {\em Ann. Sci. \'Ecole Norm. Sup. (4) 12}, 4 (1979), 501--661.

\bibitem{Jannsen1988}
{\sc Jannsen, U.}
\newblock Continuous \'etale cohomology.
\newblock {\em Math. Ann. 280}, 2 (1988), 207--245.

\bibitem{Katz1979}
{\sc Katz, N.~M.}
\newblock Slope filtration of {$F$}-crystals.
\newblock In {\em Journ\'ees de {G}\'eom\'etrie {A}lg\'ebrique de {R}ennes
  ({R}ennes, 1978), {V}ol. {I}}, vol.~63 of {\em Ast\'erisque}. Soc. Math.
  France, Paris, 1979, pp.~113--163.

\bibitem{KatzMessing1974}
{\sc Katz, N.~M., and Messing, W.}
\newblock Some consequences of the {R}iemann hypothesis for varieties over
  finite fields.
\newblock {\em Invent. Math. 23\/} (1974), 73--77.

\bibitem{Kedlaya2004}
{\sc {Kedlaya}, K.~S.}
\newblock {Full faithfulness for overconvergent $F$-isocrystals.}
\newblock In {\em {Geometric aspects of Dwork theory. Vol. I, II}}. Berlin:
  Walter de Gruyter, 2004, pp.~819--835.

\bibitem{TrihanMatsuda2004}
{\sc {Matsuda}, S., and {Trihan}, F.}
\newblock {Image directe sup\'erieure et unipotence.}
\newblock {\em {J. Reine Angew. Math.} 569\/} (2004), 47--54.

\bibitem{McCarthy1997}
{\sc McCarthy, R.}
\newblock Relative algebraic {$K$}-theory and topological cyclic homology.
\newblock {\em Acta Math. 179}, 2 (1997), 197--222.

\bibitem{Milne1980}
{\sc Milne, J.~S.}
\newblock {\em \'{E}tale cohomology}, vol.~33 of {\em Princeton Mathematical
  Series}.
\newblock Princeton University Press, Princeton, N.J., 1980.

\bibitem{Morrow_Deformational_Hodge}
{\sc Morrow, M.}
\newblock A case of the deformational {H}odge conjecture via a pro
  {H}ochschild--{K}ostant--{R}osenberg theorem.
\newblock {\em Comptes Rendus Math\'ematiques 352}, 3 (2013), 173--177.

\bibitem{Ogus1984}
{\sc {Ogus}, A.}
\newblock {F-isocrystals and de Rham cohomology. II: Convergent isocrystals.}
\newblock {\em {Duke Math. J.} 51\/} (1984), 765--850.

\bibitem{Ogus1990}
{\sc Ogus, A.}
\newblock The convergent topos in characteristic {$p$}.
\newblock In {\em The {G}rothendieck {F}estschrift, {V}ol.\ {III}}, vol.~88 of
  {\em Progr. Math.} Birkh\"auser Boston, Boston, MA, 1990, pp.~133--162.

\bibitem{Petrequin2003}
{\sc {Petrequin}, D.}
\newblock {Chern classes and cycle classes in rigid cohomology. (Classes de
  Chern et classes de cycles en cohomologie rigide).}
\newblock {\em {Bull. Soc. Math. Fr.} 131}, 1 (2003), 59--121.

\bibitem{Popescu1985}
{\sc Popescu, D.}
\newblock General {N}\'eron desingularization.
\newblock {\em Nagoya Math. J. 100\/} (1985), 97--126.

\bibitem{Popescu1986}
{\sc Popescu, D.}
\newblock General {N}\'eron desingularization and approximation.
\newblock {\em Nagoya Math. J. 104\/} (1986), 85--115.

\bibitem{Shiho2002}
{\sc {Shiho}, A.}
\newblock {Crystalline fundamental groups. II: Log convergent cohomology and
  rigid cohomology.}
\newblock {\em {J. Math. Sci., Tokyo} 9}, 1 (2002), 1--163.

\bibitem{Shiho2007}
{\sc Shiho, A.}
\newblock On logarithmic {H}odge-{W}itt cohomology of regular schemes.
\newblock {\em J. Math. Sci. Univ. Tokyo 14}, 4 (2007), 567--635.

\bibitem{Shiho2007I}
{\sc Shiho, A.}
\newblock Relative log convergent cohomology and relative rigid cohomology {I}.
\newblock {\em {\tt arXiv:0707.1742}\/} (2007).

\bibitem{Shiho2007II}
{\sc Shiho, A.}
\newblock Relative log convergent cohomology and relative rigid cohomology
  {II}.
\newblock {\em {\tt arXiv:0707.1743}\/} (2007).

\bibitem{Shiho2008}
{\sc Shiho, A.}
\newblock Relative log convergent cohomology and relative rigid cohomology
  {III}.
\newblock {\em {\tt arXiv:0805.3229}\/} (2007).

\bibitem{Tate1994}
{\sc {Tate}, J.}
\newblock {Conjectures on algebraic cycles in $\ell$-adic cohomology.}
\newblock In {\em {Motives. Proceedings of the summer research conference on
  motives, held at the University of Washington, Seattle, WA, USA, July
  20-August 2, 1991}}. Providence, RI: American Mathematical Society, 1994,
  pp.~71--83.

\end{thebibliography}

\def\cprime{$'$}

\noindent Matthew Morrow\hfill{\tt morrow@math.uni-bonn.de} \\
Mathematisches Institut\hfill \url{http://www.math.uni-bonn.de/people/morrow/}\\
Universit\"at Bonn\\
Endenicher Allee 60\\
53115 Bonn, Germany
\end{document}